\documentclass{m2an}

\usepackage[colorlinks,linkcolor=blue,anchorcolor=blue,citecolor=red]{hyperref}

\newtheorem{theorem}{Theorem}

\newtheorem{lemma}[theorem]{Lemma}
\newtheorem{remark}{Remark}

\usepackage{bm}
\usepackage{float}
\usepackage{multirow}
\usepackage{upgreek}
\usepackage{amsfonts}
\usepackage{amssymb}
\usepackage{graphicx}
\usepackage{mathrsfs}
\usepackage{tabularx}
\usepackage{array}
\usepackage{cite}
\usepackage{amsthm}
\usepackage{amsmath}

\begin{document}
\title{Hypocoercivity for the linear semiconductor 
Boltzmann equation with boundaries and uncertainties}
\author{Hongxu Chen}\address{The Chinese University of Hong Kong, Hong Kong}
\author{Liu Liu}\address{The Chinese University of Hong Kong, Hong Kong}
\author{Jiayu Wan}\address{The Chinese University of Hong Kong, Hong Kong}

\begin{abstract}
In this paper, we establish hypocoercivity for the semiconductor Boltzmann equation with the presence of an external electrical potential under the Maxwell boundary condition. We will construct a modified entropy Lyapunov functional, which is proved to be equivalent to some weighted norm of the corresponding function space. We then show that the entropy functional dissipates along the solutions, and the exponential decay to the equilibrium state of the system follows by a Gronwall-type inequality. We also generalize our arguments to situations where uncertainties in our model arise, and the hypocoercivity method we have established is adopted to analyze the regularity of the solutions along the random space.  
\end{abstract}

\date{}

\maketitle

\section{Introduction}

Kinetic equations have been widely used in many areas such as mechanics, rarefied gas, plasma physics, astrophysics, semiconductor device modeling, and social and biological sciences \cite{Villani02}. They describe the non-equilibrium dynamics of a system composed of a large number of particles and bridge atomistic and continuum models in the hierarchy of multiscale modeling. The Boltzmann-type equation, as one of the most representative models in kinetic theory, provides a powerful tool to describe molecular gas dynamics, radiative transfer, plasma physics, and polymer flow \cite{A15}. In particular, the linear semiconductor Boltzmann equation we will study in this article describes the evolution of electron density(in the phase space) within semiconductor devices, in the presence of an external electrical potential and scatterings among various particles \cite{Jungel}. 

A fundamental problem in the asymptotic analysis of kinetic equations is to show the exponential decay to the equilibrium state of the system. Among the methods proposed in this direction, \textit{hypocoercivity} methods have garnered increasing attention in recent years. The key of these methods is to build a modified entropy Lyapunov functional $H(f)$, that is equivalent to the specific norm under consideration, and to demonstrate that this functional dissipates along the solutions. From there, an explicit rate of convergence to the equilibrium state can be computed, which is one of the main advantages of the methods. A comprehensive introduction to hypocoercivity methods can be found in \cite{herau2018}. For various theoretical results related to hypocoercivity, we refer the reader to \cite{villani2009, MN2006, briant2015, herau2006, HN2004}. Most of these results establish exponential convergence in some weighted $H^1$ space and require certain regularity on the initial datum. Recently, Dolbeault, Mouhot and Schmeiser proposed a hypocoercivity method in \cite{DM2015}, which leads to exponential decay to the equilibrium state in some weighted $L^2$ space, and hence avoiding the burden of computing derivatives when generalized to the numerical settings(\cite{BHR2020}). In this work, we will focus on the framework in \cite{DM2015} and adapt it to establish hypocoercivity in our setting. 

All the hypocoercivity methods mentioned in the previous paragraph consider either a periodic spatial domain $\mathbb{T}^3$ or an unbounded domain $\mathbb{R}^3$, which fails to capture the phenomena with a physical boundary. In particular, in semiconductor devices, electrons are assumed to move within a bounded container and interact in a certain way with the wall of the container. If there is no electrical charge exchange at the boundary, e.g., the wall is insulating, the electrons are simply reflected at the wall, with possible changes in momentum and energy while maintaining the total charges. In the context of the Boltzmann equation for rarefied gas, Guo established the $L^2$ hypocoercivity for systems with inflow, diffuse reflection, and specular reflection boundary conditions using a non-constructive method \cite{Guo2010}. Later developments, including \cite{esposito2013non, BCMT2023, esposito2018stationary, briant2016asymptotic}, introduced constructive approaches to hypocoercivity, which have been applied to study a wide class of kinetic equations with physical boundaries (see, e.g., \cite{CKP2025, guo2020landau}). These works also address boundary value problems involving the Maxwell boundary condition, which serves as an intermediate model between diffuse and specular reflection. For more literature on the Maxwell boundary and generalized diffuse boundary conditions, we refer to \cite{mischler2010, chen, chen2025boltzmann}. In this paper, we establish hypocoercivity for the semiconductor Boltzmann equation under the Maxwell boundary condition by constructing a modified entropy Lyapunov functional.

Another important physical phenomenon that we need to address in the analysis of kinetic equations is uncertainty quantification. Kinetic equations typically involve various sources of uncertainty, such as modeling errors, imprecise measurements, and uncertain initial conditions. In particular, the scattering kernel in the linear semiconductor Boltzmann equation governs the transition rates during particle collisions. Calculating this scattering kernel from first principles is highly complex, and in practice, heuristic approximations or empirical data are often used, inevitably introducing uncertainties. As a result, addressing uncertainty quantification(UQ) becomes essential for evaluating, validating, and improving the underlying models. To numerically simulate kinetic models with randomness, the generalized polynomial chaos (gPC)-based stochastic Galerkin (SG) method and its variations have been widely adopted, demonstrating success in a range of applications \cite{Xiu}. For numerical experiments and theoretical analysis of gPC-based methods, we refer readers to works such as \cite{JinPareschi,parUQ, LiuJin2018, LiuDaus2019, WL2025, HJ16, LQ2024-1, LQ2024-2}. As mentioned in these works, the spectral convergence of gPC methods heavily relies on the regularity of solutions along the random space, and a general framework for analyzing this regularity is proposed in \cite{LW2017}. In this paper, we will also generalize the established hypocoercivity method to study the regularity along the random space, when the linear semiconductor model and the Maxwell boundary condition are considered.  

For the rest of the paper, we will establish hypocoercivity for linear semiconductor Boltzmann equation with the Maxwell boundary condition, based on the construction of a modified entropy Lyapunov functional. In the presence of a physical boundary, contributions from the boundary terms need to be carefully analyzed before fitting into the general framework. In addition, we will generalize the arguments to study the regularity along the random space, when uncertainties of the model are present. Our main results can be found in Theorem \ref{exponential decay of f base case}, Theorem \ref{exponential decay of f with potential} and Theorem \ref{exponential decay of g_l}. 

A short outline of the paper is given as follows. In Section 2, we present a short introduction to the semiconductor Boltzmann equation and the Maxwell boundary condition. In Section 3, we establish hypocoercivity for a base model where external electrical potential is negligible, to demonstrate the core methodology along with some critical aspects of our proof. We then establish hypocoercivity in the presence of external potentials in Section 4. Finally in Section 5, we generalize the arguments to study the regularity of solutions along the random space when the model involves some uncertainties, and some concluding remarks are provided in Section 6.

\newpage
\section{preliminaries}

In this section, we present a short introduction to linear semiconductor Boltzmann equation and summarize some key facts that will be useful for our analysis later. A discussion of the Maxwell boundary condition will also be included. To begin with, the general Boltzmann equation for semiconductor devices is given by 

\begin{equation}\label{semiconductor Boltzmann equation}
    \begin{cases}
         \partial_{t}f + v\cdot \nabla_{x}f + \nabla_{x}V \cdot \nabla_{v}f= Q(f) \\
         \ f(0,x,v) = f_{0}(x,v), \hspace{0.7cm} (x,v) \in \Omega \times \mathbb{R}^3,
    \end{cases}
\end{equation}

\noindent where $f=f(t,x,v)$ is the probability density distribution for particles(typically electrons) at $x \in \Omega$ with velocity $v \in \mathbb{R}^3$. We assume the space domain $\Omega \subset \mathbb{R}^3$ is bounded with proper boundary conditions(will be described later). $V = V(t,x)$ is the external electrical potential which is prescribed at the beginning. The collision operator $Q$ describes the scattering of particles. In semiconductor crystals, there are three main sources for scattering of electrons: phonons(electron-phonon scattering), lattice defects(ionized impurity scattering), and other carriers(carrier-carrier scattering). 

To derive a specific expression for $Q$, let $P(t,x,v_{*} \rightarrow v)$ be the rate at which a particle at $(t,x)$ changes its velocity $v_{*}$ into another velocity $v$ due to a scattering event. By Pauli exclusion principle(valid for fermions and hence for electrons), we assume $P(t,x,v_{*} \rightarrow v)$ is proportional to the occupation probability $f(t,x,v_{*})$ and the probability $1-f(t,x,v)$ that the phase state $(x,v)$ is not occupied. Hence 

\begin{equation}\label{expression for P}
    P(t,x,v_{*} \rightarrow v) = s(x,v_{*},v)f(t,x,v_{*})(1-f(t,x,v)),
\end{equation}

\noindent where $s(x,v_{*},v)$ is called the scattering rate, and we assume it is homogeneous in space: $s(x,v_{*},v)=s(v_{*},v)$. Then the rate of change of $f$ due to collision(scattering) is the sum of all in-scattering rates from some $v_{*}$ to $v$ minus the out-scattering rates from $v$ to some $v_{*}$, and hence

\begin{equation}\label{expression for Q}
    \begin{aligned}
        Q(f)(t,x,v) =& \int_{\mathbb{R}^3} P(t,x,v_{*} \rightarrow v) - P(t,x,v \rightarrow v_{*}) ~dv_{*} \\
        =& \int_{\mathbb{R}^3} s(v_{*},v)f_{*}(1-f) -s(v,v_{*})f(1-f_{*}) ~dv_{*},
    \end{aligned}
\end{equation}

\noindent where $f=f(t,x,v)$ and $f_{*}=f(t,x,v_{*})$. We further assume the principle of detailed balance: 

\begin{equation}\label{principle of detailed balance}
    \frac{s(v_{*},v)}{\mathcal{M}(v)}=\frac{s(v,v_{*})}{\mathcal{M}(v_{*})},
\end{equation}

\noindent where $\mathcal{M}(v)$ is the normalized Maxwellian given by 

\begin{equation}\label{def of normalized Maxwellian}
    \mathcal{M}(v)=\frac{1}{(2\pi)^{\frac{3}{2}}} e^{-\frac{|v|^2}{2}}.
\end{equation}

Typically, the collision operator $Q$ is nonlinear and it is hard to tackle directly. To simplify the situation, we use the low-density approximation. We assume that the distribution function $f$ is small in the sense that $ 0 \leq f(t,x,v) \ll 1$ and hence $1-f(t,x,v) \approx 1$. Then using \eqref{principle of detailed balance}, we have 

\begin{equation*}
    \begin{aligned}
        Q(f)(v) =& \int_{\mathbb{R}^3} s(v_{*},v)f_{*} - s(v,v_{*})f ~dv_{*} \\
        =& \int_{\mathbb{R}^3} \frac{s(v_{*},v)}{\mathcal{M}}(\mathcal{M}f_{*} - \mathcal{M}_{*}f) ~dv_{*}.
    \end{aligned}
\end{equation*}

\noindent We let $\sigma(v_{*},v) = \frac{s(v_{*},v)}{\mathcal{M}}$, which is called the collision cross-section, and it is symmetric by \eqref{principle of detailed balance}. Clearly, under the low-density approximation, the collision operator $Q$ becomes linear and we replace the notation $Q$ by $\mathcal{L}$ to indicate this linearity. Then we derive the linear semiconductor Boltzmann equation  

\begin{equation}\label{linear semiconductor Boltzmann equation}
    \begin{cases}
         \partial_{t}f + v\cdot \nabla_{x}f+ \nabla_{x}V \cdot \nabla_{v}f = \mathcal{L}(f) \\
         \ f(0,x,v) = f_{0}(x,v), \hspace{0.7cm} (x,v) \in \Omega \times \mathbb{R}^3,
    \end{cases}
\end{equation}

\noindent with the linear collision(scattering) operator given by 

\begin{equation}\label{expression for L}
    \mathcal{L}(f)(t,x,v) = \int_{\mathbb{R}^3}\sigma(v,v_{*})(\mathcal{M}f_{*} - \mathcal{M}_{*}f) ~dv_{*}.
\end{equation}

Since we assume the cross-section $\sigma$ is spatially independent, the operator $\mathcal{L}$ only acts on the velocity variable $v$ and is local in $t$ and $x$. It is well-known that the kernel of $\mathcal{L}$ is spanned by the normalized Maxwellian, i.e.:

\begin{equation}\label{kernel of L}
    N(\mathcal{L}) = span \{ \mathcal{M}(v)\}.
\end{equation}

\noindent In addition, it can be shown that $\mathcal{L}$ satisfies the local conservation of mass, i.e.:  

\begin{equation}\label{local conservation of mass}
    \int_{\mathbb{R}^3} \mathcal{L}(f) ~dv =0.
\end{equation}

\noindent We remark that typically $\mathcal{L}$ does not conserve momentum and energy(for example, phonon emission or absorption during the scattering will result in a loss or gain of total energy). We therefore define 

\begin{equation}\label{def of j_L}
    j^{\mathcal{L}}:= \int_{\mathbb{R}^3} v\mathcal{L}(f) ~dv,
\end{equation}

\noindent which is an important quantity for our analysis later. 

The precise expression for the cross-section $\sigma(v,v_{*})$ depends on the exact underlying scattering event, and it is usually quite complicated. Fortunately, it is not necessary for us to know the details of the scattering mechanism, and some assumptions about $\sigma(v,v_{*})$ will suffice for our purpose. We thus make the following two assumptions on the collision cross-section $\sigma(v,v_{*})$: 

\bigskip

\noindent \textbf{Assumption 1}. 
There exists $\lambda>0$ such that
\begin{equation}\label{H1}
    \sigma(v,v_{*}) \geq \lambda >0.
\end{equation}

\noindent \textbf{Assumption 2}. There exists $C<\infty$ such that

\begin{equation}\label{H2}
    \iint_{\mathbb{R}^3 \times \mathbb{R}^3} \sigma(v,v_{*})^{2}\mathcal{M}\mathcal{M}_{*}(|v|^2 + |v_{*}|^2) ~dvdv_{*} \leq C .
\end{equation}

Based on these two assumptions, we can derive some important estimates regarding the collision operator $\mathcal{L}$. Let $\langle \cdot, \cdot \rangle_{dx}$ be the standard inner product on $L^2(\Omega)$, and define $\langle \cdot , \cdot \rangle_{d\nu}$ to be the weighted inner product on 
\begin{align*}
    L^2(d\nu):=\left\{ f \Big| \int_{\mathbb{R}^3} f^2 ~d\nu < \infty \right\}, \text{ with the measure } d\nu := \mathcal{M}^{-1}dv .
\end{align*}
 Let $\| \cdot \|_{dx}~, ~\| \cdot \|_{d\nu}$ be the associated norms over the corresponding function spaces. Define $\pi_{\mathcal{L}}$ to be the projection operator with respect to $\langle \cdot , \cdot \rangle_{d\nu}$ onto the space spanned by the unit vector $\mathcal{M}(v)$ in $L^{2}(d\nu)$, ie:

\begin{equation}\label{def of pi_L}
    \pi_{\mathcal{L}}(f) = \left(\int_{\mathbb{R}^3} f ~dv\right)\mathcal{M}(v) = \rho(t,x) \mathcal{M}(v).
\end{equation}

\noindent where $\rho(t,x)= \int_{\mathbb{R}^3} f ~dv$ is the mass density of particles. We let $f^{\perp}:=f - \pi_{\mathcal{L}}(f)$. Then, we have the following two important estimates, whose proofs are included in the appendix: 

\begin{lemma}\label{microscopic coercivity}
    \textnormal{(microscopic coercivity).} Suppose \textbf{Assumption 1} holds, then 

    \begin{equation}\label{microscopic coercivity inequality}
        \langle \mathcal{L}(f) , f \rangle_{d\nu} \leq -\lambda \| f^{\perp }\|^2_{d\nu}.
    \end{equation}
\end{lemma}

\begin{lemma}\label{estimate for j_L}
    \textnormal{(estimate for $j^\mathcal{L}$).} Let $j^{\mathcal{L}}$ be defined as in \eqref{def of j_L}. Suppose \textbf{Assumption 2} holds, then for some $C_{\mathcal{L}}>0$, it holds that 

    \begin{equation}\label{j_L inequality}
        \| j^{\mathcal{L}} \|_{dx} \leq C_{\mathcal{L}} \| f^{\perp} \|_{dxd\nu}.
    \end{equation}

\noindent 
\end{lemma}

Now, we describe the boundary condition for the problem. We denote the boundary of the phase space as 

\begin{equation*}
    \gamma := \{ (x,v) \in \partial \Omega \times \mathbb{R}^3\}. 
\end{equation*}

\noindent Let $n(x)$ be the outward normal direction at $x \in \partial \Omega$, we decompose $\gamma$ as 

\begin{equation}\label{boundary decomposition}
    \begin{aligned}
        \gamma_{-}=& \{ (x, v) \in \partial \Omega \times \mathbb{R}^3: n(x) \cdot v < 0 \}, \\
        \gamma_{+}=& \{ (x, v) \in \partial \Omega \times \mathbb{R}^3: n(x) \cdot v > 0 \}, \\
        \gamma_{0}=& \{ (x, v) \in \partial \Omega \times \mathbb{R}^3: n(x) \cdot v = 0 \}.
    \end{aligned}
\end{equation}

\noindent Maxwell boundary condition is prescribed for the incoming phase space, given by
\begin{align}
    &   f(t,x,v)|_{\gamma_-} = c f(t,x,v-2(n(x)\cdot v)n(x))  + (1-c) C_M \mathcal{M}(v)\int_{n(x)\cdot u>0} f(t,x,u)(n(x)\cdot u) du,  \label{maxwell_bc}
\end{align}
where $0\leq c\leq 1$, $C_M=\sqrt{2\pi}$ is chosen so that $C_M\mathcal{M}(v)|n(x) \cdot v|$ is a probability measure on the half velocity space $\{v \in \mathbb{R}^3: n(x) \cdot v <0 \}$, i.e.

\begin{equation}\label{unity of c_mu}
    \int_{n(x) \cdot v <0} C_M \mathcal{M}(v) |n(x) \cdot v| ~dv = 1.
\end{equation}

When $c=1$, it corresponds to the specular boundary condition; when $c=0$, it corresponds to the diffusive boundary condition. Physically, under the specular reflection, the distribution of the re-emitted particles is fully determined by that of the incoming particles. While under the diffusive reflection, the re-emitted particles have completely lost memory of the incoming ones. Maxwell boundary condition is simply a convex combination of these two reflections, and it typically occurs when the spatial boundary $\partial \Omega$ is insulating in the sense that no carrier can enter or exit the boundary but they are simply reflected. This boundary condition exhibits the null-flux property, which is stated as follows and the proof is provided in the Appendix. 

\begin{lemma}\label{null-flux property}
    \textnormal{(null-flux property).} If $f$ satisfies the Maxwell boundary condition, then 

    \begin{equation}\label{null-flux}
        \int_{\mathbb{R}^3} f(t, x, v) \left( n(x) \cdot v \right) ~dv = 0 .
    \end{equation}
\end{lemma}

Finally, in this section, we introduce a technical tool to deal with boundary functions, which will be crucial for our proofs later. For fixed $x \in \partial \Omega$, we define the boundary inner product over $\gamma_{+}$ in $v$ as:

\begin{equation}\label{inner product on boundary}
    \langle g_1, g_{2} \rangle_{\gamma_{+}}=\int_{n(x) \cdot v > 0} g_1 g_2  (n(x) \cdot v) ~dv,
\end{equation}

\noindent for $g_1, g_2$ defined over $\partial \Omega \times \mathbb{R}^3$. we define an ${L_v^2}$ - projection for a boundary function $h(x, v)$ towards $\sqrt{C_M}\sqrt{\mathcal{M}(v)}$ with respect to $\langle , \rangle_{\gamma_{+}}$ as:

\begin{equation}\label{projection on boundary}
    P_{\gamma} h=C_M\sqrt{ \mathcal{M}(v)}\int_{n(x) \cdot u>0} h  \sqrt{ \mathcal{M}(u)}  (n(x) \cdot u)~ du  .
\end{equation}

\bigskip

\section{hypocoercivity for the base case}
In this section, we start to use hypocoercivity methods to prove the exponential decay of the solution $f$ to \eqref{linear semiconductor Boltzmann equation} in some properly normed functional space. We begin with the most basic case where there is no uncertainty in our system and the external electrical potential $V$ can be ignored. We then generalize our arguments include uncertainties and electrical potential in the next two sections. So the Boltzmann equation we will be considering in this section is given by 

\begin{equation}\label{linear semiconductor base case}
    \begin{cases}
         \partial_{t}f + v\cdot \nabla_{x}f = \mathcal{L}(f) \\
         \ f(0,x,v) = f_{0}(x,v), \hspace{0.7cm} (x,v) \in \Omega \times \mathbb{R}^3,
    \end{cases}
\end{equation}

\noindent where $\mathcal{L}$ is given by \eqref{expression for L} such that the cross-section $\sigma(v,v_{*})$ satisfies \textbf{Assumption 1} and \textbf{Assumption 2}, and we assume Maxwell boundary condition on $\partial \Omega \times \mathbb{R}^3$ with Lipschitz boundary $\partial \Omega$. Our main result of this section is the following Theorem. 

\begin{theorem}\label{exponential decay of f base case}
    Let $f$ be a solution to \eqref{linear semiconductor base case}, where the cross-section $\sigma(v,v_{*})$ of $\mathcal{L}$ satisfies \textbf{Assumption 1} and \textbf{Assumption 2}. Suppose $f$ satisfies the Maxwell boundary condition \eqref{maxwell_bc} with parameter $c$ on $\gamma_-$. Then there exists $C>0$ and $\tau >0$ such that for $\forall t \geq 0$, 

    \begin{equation}\label{exponential decay of f equation base case}
        \| f(t) - \mu_{f}\mathcal{M}\|_{dxd\nu} \leq C \| f_{0} -\mu_{f}\mathcal{M} \|_{dxd\nu} ~e^{-\tau t},
    \end{equation}

\noindent where $\mu_{f} = \iint_{\Omega \times \mathbb{R}^3} f_{0} ~dxdv$.  
\end{theorem}

Notice that $\mathcal{L}$ satisfies the local conservation of mass (\eqref{local conservation of mass}), and there is no mass flow along the boundary by \eqref{null-flux}. So the total mass of the system is conserved, ie: 

\begin{equation}\label{global conservation of mass}
    \frac{d}{dt} \iint_{\Omega \times \mathbb{R}^3} f(t,x,v) ~dxdv =0,
\end{equation}

\noindent and hence $\mu_{f} = \iint_{\Omega \times \mathbb{R}^3} f_{0} ~dxdv = \iint_{\Omega \times \mathbb{R}^3} f(t) ~dxdv ~~\forall t \geq 0$. This is why we have the same $\mu_{f}$ on both sides of \eqref{exponential decay of f equation base case}. Also notice that adding $a\mathcal{M}(v)$ for any $a \in \mathbb{R}$ to $f$ on both sides of \eqref{linear semiconductor base case} will not change the equation due to the linearity of $\mathcal{L}$ and \eqref{kernel of L}, and there is no change on the boundary condition as well since $\mathcal{M}(v)$ satisfies the Maxwell boundary condition. So we can replace $f$ by $f+a\mathcal{M}$ for some properly chosen $a$ so that $\mu_{f}=0$. For the rest of the section, we assume $\mu_{f}=0$ to simplify the arguments.  

To prove the exponential decay of $f$, we borrow the technique introduced by \cite{DM2015}, which is reformulated in \cite{BHR2020} in terms of moments method. We will then generalize the arguments to fit in the Maxwell boundary condition we are interested in(both \cite{DM2015} and \cite{BHR2020} only consider $\mathbb{T}^3$ or $\mathbb{R}^3$ as the spatial domain). The idea is to use the moments of $f$ to construct a modified entropy functional $H(f)$, which is equivalent to $\| f \|_{dxd\nu}^2$. Then apply various estimates to show 

\begin{equation*}
    \frac{d}{dt}H(f) \leq -\tau H(f),
\end{equation*}

\noindent and hence the result follows by a Gronwall-type inequality. For the rest of the section, we will discuss the construction of $H(f)$ and the estimates in detail, and we conclude the section with a proof of Theorem \ref{exponential decay of f base case}. To simplify the notations, we use $\langle \cdot, \cdot \rangle$ to stand for $\langle \cdot, \cdot \rangle_{dxd\nu}$ and $\| \cdot \|$ to stand for $\| \cdot \|_{dxd\nu}$. 

We start with moments and the moment equation of $f$. Define the moments of $f$ up to order 2 as
\begin{equation}\label{moments of f}
    \begin{cases}
        & \rho(t,x) = \int_{\mathbb{R}^3} f ~dv \\
        & j(t,x) = \int_{\mathbb{R}^3} vf ~dv \\
        & S(t,x) = \int_{\mathbb{R}^3} (v \otimes v - I)f(t,x,v) ~dv, \hspace{0.7cm} v \otimes v := vv^{T} = \left(\begin{array}{lll}
                            v_1v_1 & v_1v_2 & v_1v_3 \\
                            v_2v_1 & v_2v_2 & v_2v_3 \\
                            v_3v_1 & v_3v_2 & v_3v_3 \\ 
                            \end{array}\right).
    \end{cases}
\end{equation}

\noindent  Multiplying equation \eqref{linear semiconductor base case} by 1, $v$ and taking integration, we get

\begin{equation}\label{moment equations}
   \begin{cases}
        \partial_t \rho + \nabla_x \cdot j=0 \\
        \partial_t j+\nabla_x \cdot(S + \rho I) = j^{\mathcal{L}},
    \end{cases}
\end{equation}

\noindent where $j^{\mathcal{L}}=\int_{\mathbb{R}^3} v \mathcal{L}(f) ~dv$ is as defined in \eqref{def of j_L} and $I$ is the identity matrix in $M_{3}(\mathbb{R})$. We use the notation $\nabla_x \cdot A$ to stand for $\left(\begin{array}{l} \nabla_x \cdot A^1  \\ \nabla_x \cdot A^2  \\ \nabla_x \cdot A^3  \\ \end{array}\right)$ for a matrix $A$ and $A^i$ refers to the ith row of $A$. Let $\tilde{S}:=\int_{\mathbb{R}^3} (v \otimes v) f ~dv$ so that $S=\tilde{S} - \rho I$. We have the following estimates for the moments of $f$: 

\begin{lemma}\label{moment estimates}
    \textnormal{(moment estimates).} 
For some $D>0$, it holds that
    \begin{equation*}
        \begin{aligned}
            & (1).~ \| \rho \|_{dx} \leq \| f \| ,\\
            & (2).~ \| j \|_{dx}  \leq \sqrt{3}\| f^{\perp} \| ,\\
            & (3).~ \| S \|_{dx} \leq D \| f^{\perp} \|.
        \end{aligned}
    \end{equation*}

\end{lemma}

\textit{Proof.} For the first estimate, we have

\begin{equation*}
    \| \rho \|_{dx} = \| \langle f, \mathcal{M}\rangle_{d\nu} \|_{dx} \leq \| \|f\|_{d\nu} \|_{dx} = \| f \|
\end{equation*}

\noindent using Cauchy-Schwarz and the fact that $\| \mathcal{M}(v)\|_{d\nu}=1$. 

For the second estimate, we have

\begin{equation*}
    j=\int_{\mathbb{R}^3} v f ~dv=\left(\begin{array}{l}\int_{\mathbb{R}^3} v_1 f ~dv \\ \int_{\mathbb{R}^3} v_2 f ~dv \\ \int_{\mathbb{R}^3} v_3 f ~dv\end{array}\right) = \left(\begin{array}{l}\int_{\mathbb{R}^3} v_1 (\rho\mathcal{M} + f^{\perp}) ~dv \\ \int_{\mathbb{R}^3} v_2 (\rho\mathcal{M} + f^{\perp}) ~dv \\ \int_{\mathbb{R}^3} v_3 (\rho\mathcal{M} + f^{\perp}) ~dv\end{array}\right) = \left(\begin{array}{l}\int_{\mathbb{R}^3} v_1 f^{\perp} ~dv \\ \int_{\mathbb{R}^3} v_2 f^{\perp} ~dv \\ \int_{\mathbb{R}^3} v_3 f^{\perp} ~dv\end{array}\right) = \left(\begin{array}{l} \langle f^{\perp}, v_1\mathcal{M}\rangle_{d\nu} \\ \langle f^{\perp}, v_1\mathcal{M}\rangle_{d\nu} \\ \langle f^{\perp}, v_1\mathcal{M}\rangle_{d\nu} \end{array}\right)
\end{equation*}

\noindent using the fact that $\int_{\mathbb{R}^3} v_i \mathcal{M}(v) ~dv=0$ for $1 \leq i \leq 3$. Then using Cauchy-Schwarz and the fact that $\| v_{i}\mathcal{M}\|=1$ for $1 \leq i \leq 3$, we have $\langle f^{\perp}, v_{i}\mathcal{M}\rangle_{d\nu} \leq \| f^{\perp} \|_{d\nu}$ for $1 \leq i \leq 3$, and hence $\| j\|_{dx} \leq \sqrt{3}\| f^{\perp} \|$.

For the estimate of $S$, we notice that 

\begin{equation*}
    S  =\int_{\mathbb{R}^3}(v \otimes v-I) f ~dv = \int_{\mathbb{R}^3}(v \otimes v-I) (\rho \mathcal{M} + f^{\perp}) ~dv = \int_{\mathbb{R}^3}(v \otimes v-I) f^{\perp} ~dv
\end{equation*}

\noindent using the fact that $\int_{\mathbb{R}^3} v_i v_j \mathcal{M}(v) ~dv=\delta_{i j}$, where $\delta_{i j}$ is the Kronecker delta symbol for $ 1 \leq i,j \leq 3$. Then we have $\int_{\mathbb{R}^3}(v \otimes v-I) f^{\perp} ~dv = \langle(v \otimes v-I) \mathcal{M}, f^{\perp}\rangle_{d\nu}$. Let $D = \| (v \otimes v-I) \mathcal{M} \|_{d\nu}$, then $\| S \|_{dx} \leq D \| f^{\perp} \|$ by Cauchy-Schwarz inequality. 

\qed

Now, we construct the ``modified entropy functional" $H(f)$, which will serve as the key ingredient of our analysis. To this end, we let 

\begin{equation}\label{def of H(f)}
    H(f)=\frac{1}{2}\| f \|^2 + \eta\langle j, \nabla_x \phi \rangle_{d x},
\end{equation}

\noindent where $\eta >0$ is a small parameter to be chosen and $\phi(t,x)$ is the solution of the Poisson equation with Neumann boundary condition:
\begin{equation}\label{Poisson equation for phi}
\begin{cases}
   & -\Delta \phi =\rho, \hspace{0.7cm} \int_{\Omega} \phi ~dx =0 \\
   & \nabla_{x}\phi \cdot n(x)=0  \ \text{ for } \ \forall x \in \partial \Omega.
\end{cases}    
\end{equation}

\noindent We comment that we require $\phi$ to satisfy the Neumann boundary condition so that it is compatible with the Maxwell boundary condition of $f$, as we will see later. Note that the existence of the solution to the Poisson equation is established by our assumption that $\int_{\Omega}\rho(x,t) ~dx =0 ~\forall t$, and the uniqueness of the solution is guaranteed by the additional requirement $\int_{\Omega} \phi ~dx=0$. Since we assume the boundary $\partial \Omega$ is Lipschitz, Poisson regularity estimate and trace theorem hold for $\phi$, i.e, for some $K, D_{\gamma}>0$,  

\begin{equation*}
    \begin{cases}
        \| \phi\|_{H^2(\Omega)} \leq K\| \rho \|_{L^2(\Omega)}=K\| \rho \|_{dx} \\
        \Vert\nabla_{x}\phi \Vert_{L^2(\partial \Omega)} \leq D_{\gamma}\| \nabla_{x}\phi\|_{H^1(\Omega)}.
    \end{cases}
\end{equation*}

\noindent 

\begin{remark}
    \eqref{Poisson equation for phi} is a slight simplification of the original version presented in \cite{DM2015}, where $\phi$ satisfies a more general elliptic equation instead of a Poisson equation. We will stick to \eqref{Poisson equation for phi} since it is simpler and it suffices for our purpose. 
\end{remark}

We now derive some estimates for $\phi$, which will be important for our analysis: 

\begin{lemma}\label{estimates_for_phi}
    \textnormal{(Estimates for $\phi$).} For $\phi$ described above, we have the following 
    
\begin{equation*}
    \begin{aligned}
        &\text{(1). } \left\|\nabla_x \phi\right\|_{d x} \leq C_p \| \rho \|_{dx} ,\\
        &\text{(2). } \| \partial_{t} \nabla_{x}\phi \|_{dx} \leq \| j \|_{dx}.
    \end{aligned}
\end{equation*}
\end{lemma}

\textit{Proof.} For the first estimate, we apply integration by parts to yield

\begin{equation*}
    \begin{aligned}
        \|\nabla_x \phi \|_{dx}^2 =& \langle \nabla_x \phi, \nabla_x \phi \rangle_{dx} \\
         =&- \langle \phi, \underbrace{\nabla_x \cdot \nabla_x \phi}_{=\Delta \phi = -\rho} \rangle_{dx} + \int_{\partial \Omega} \phi \underbrace{(\nabla_x \phi \cdot n(x))}_{=0 \text{  by Neumann boundary condition on } \phi} ~d S_x \\
         =& \langle \phi, \rho \rangle_{dx} 
         \leq \|\phi\|_{dx} \cdot\| \rho \|_{dx} \\
         \leq& C_{p} \| \nabla_{x} \phi\|_{dx} \cdot \| \rho \|_{dx},
    \end{aligned}
\end{equation*}

\noindent where we apply Cauchy-Schwarz in the second last inequality and the Poincaré inequality in the last line. Canceling out $\| \nabla_{x} \phi \|$ on both sides gives the first estimate. 

For the second estimate, we again apply integration by parts to get 
\begin{equation*}
    \begin{aligned}
        \|\partial_t \nabla_x \phi \|_{dx}^2=& \langle \partial_t \nabla_x \phi, \partial_t \nabla_x \phi \rangle_{dx} 
        =\langle \nabla_x \partial_t \phi, \nabla_x \partial_t \phi \rangle_{dx} \\
        =&-\langle \partial_t \phi, \underbrace{\nabla_x \cdot \nabla_x \partial_t \phi}_{=\Delta \partial_{t}\phi = \partial_{t}\Delta \phi = -\partial_{t}\rho} \rangle_{dx} + \int_{\partial \Omega} \partial_t \phi \underbrace{(\nabla_x \partial_t \phi \cdot n(x))}_{=0 \text{ by Neumann boundary condition on } \phi} ~d S_x \\
        = & \langle \partial_t \phi, \partial_t \rho \rangle_{dx} \\
        = & \langle\partial_t \phi,-\nabla_x \cdot j \rangle_{dx} \hspace{0.7cm}\text{(by the first moment equation)}\\
        = & \langle\nabla_x \partial_t \phi, j \rangle_{dx} + \int_{\partial \Omega} \partial_t \phi \underbrace{(j \cdot n(x))}_{=0 \text{ by zero-flux property of the Maxwell boundary condition}} ~d S_x \\
        = & \langle \partial_{t} \nabla_{x} \phi , j\rangle_{dx} 
        \leq \| \partial_{t} \nabla_{x} \phi \|_{dx} \cdot \| j \|_{dx},
    \end{aligned}
\end{equation*}

\noindent where we apply Cauchy-Schwarz in the last inequality. Cancelling out $\| \partial_{t} \nabla_{x} \phi \|_{dx}$ on both sides gives the second estimate. 

\qed

Now, based on Lemma \ref{moment estimates} and \ref{estimates_for_phi}, we can prove the equivalence between $H(f)$ and $\| f\|^2$:

\begin{lemma}\label{equivalent norm}
    \textnormal{(Equivalence of norms).} Let $H(f)$ be defined as \eqref{def of H(f)}, then there exists $\eta>0$ small enough, such that $H(f)$ and $\| f \|^2$ are equivlent. 
\end{lemma}

\textit{Proof}. Using Lemma \ref{moment estimates} and \ref{estimates_for_phi}, we have 

\begin{equation*}
    \eta | \langle j, \nabla_x \phi \rangle_{dx} | \leq \eta\|j\|_{dx} \cdot \|\nabla_x \phi\|_{dx} \leq \sqrt{3}\eta C_p \|f^{\perp}\| \cdot \| \rho \|_{d x} \leq \sqrt{3}\eta C_{p}(\underbrace{\| f^{\perp} \|^2 + \| \rho \mathcal{M} \|^2}_{ = \|f\|^2}), 
\end{equation*}

\noindent and hence 

\begin{equation*}
    -\sqrt{3}\eta C_p \|f\|^2 \leq \eta \langle j, \nabla_x \phi \rangle_{dx} \leq \sqrt{3}\eta C_p \|f\|^2 .
\end{equation*}
    
Choose $\eta >0$ such that $ \eta <\frac{1}{2\sqrt{3}C_p}$, then 

\begin{equation}\label{equivalent norm equation}
    c_\eta \|f\|^2 \leq H(f) \leq C_\eta \|f\|^2 ,
\end{equation}

\noindent where $c_{\eta} = \frac{1}{2}-\sqrt{3}\eta C_p >0$ and $C_{\eta} = \frac{1}{2} +\sqrt{3}\eta C_p >0$. Hence $H(f)$ and $\|f\|^2$ are equivalent. 

\qed

Before we present the proof of Theorem \ref{exponential decay of f base case}, we need one last lemma, which gives the estimate for $\frac{d}{dt}\| f \|^2$. We comment that in the proof of this lemma, one needs to be careful when performing integration by parts, since the boundary term arises from integration by parts won't necessarily vanish in the case of Maxwell boundary condition(unlike periodic boundary condition). 

\begin{lemma}\label{entropy estimate}
    \textnormal{(Entropy estimate).} Let $f$ be a solution to \eqref{linear semiconductor base case} with Maxwell boundary condition, and let $\lambda>0$ be defined as in \textbf{Assumption 1}. Then 

    \begin{equation}\label{entropy estimate equation}
        \frac{1}{2} \frac{d}{dt} \| f \|^2 \leq - \frac{1-c^2}{2}| (I-P_{\gamma}) h|_{2,+}^{2} -\lambda \| f^{\perp} \|^2 ,
    \end{equation}

    \noindent where $P_{\gamma}$ is the projection defined in \eqref{projection on boundary}, $h:=\mathcal{M}^{-\frac{1}{2}}f$, 
    \begin{align*}
        &|h|_{2,+}^{2} := \int_{\partial \Omega}\int_{n(x)\cdot v>0} |h(t,x,v)|^2 (n(x)\cdot v)~dvd S_{x}, 
    \end{align*}
$dS_{x}$ is the surface integral.
\end{lemma}

\textit{Proof.} Take $\langle, f \rangle$ on both sides of \eqref{linear semiconductor Boltzmann equation}, we get:

\begin{equation}\label{entropy equation}
    \frac{1}{2} \frac{d}{dt}\|f\|^2 + \langle v \cdot \nabla_x f, f \rangle=\langle \mathcal{L}(f), f \rangle. 
\end{equation}

For the term $\langle v \cdot \nabla_x f, f \rangle$, we have the following

\begin{equation*}
    \begin{aligned}
        \langle v \cdot \nabla_x f, f \rangle = & \int_{\mathbb{R}^3} \int_{\Omega}(\nabla_x \cdot (vf)) f \mathcal{M}^{-1}(v) ~dxdv \\
        = & - \underbrace{\int_{\mathbb{R}^3} \int_{\Omega} (vf) \cdot (\nabla_x f \mathcal{M}^{-1}(v)) ~dxdv}_{=\langle v \cdot \nabla_x f, f \rangle} + \int_{\mathbb{R}^3} \int_{\partial \Omega} f^2 \mathcal{M}^{-1}(v)(n(x) \cdot v) ~dS_{x}dv,
\end{aligned}
\end{equation*}

\noindent where we use integration by parts in the last equality. Let $f=\mathcal{M}^{\frac{1}{2}} h$, then for each $x \in \partial \Omega$, we have

\begin{equation*}
    \begin{aligned}
        \int_{\mathbb{R}^3} f^2 \mathcal{M}^{-1}(v) (n(x) \cdot v) ~dv= & \int_{\mathbb{R}^3} h^2 (n(x) \cdot v) ~dv \\
        = & \int_{n(x) \cdot v>0} h^2 (n(x) \cdot v) ~dv + \int_{n(x) \cdot v<0} h^2 (n(x) \cdot v) ~dv.
\end{aligned}
\end{equation*}

By definition \eqref{inner product on boundary}, $\int_{n(x) \cdot v >0} h^2 (n(x) \cdot v) ~dv=  \|h\|_{\gamma_{+}}^2$ where $\| \cdot\|_{\gamma_{+}}$ is the norm induced by \eqref{inner product on boundary}. Then by the Maxwell boundary condition, 

\begin{equation*}
    f(t, x, v)|_{\gamma_{-}}=cf(t,x,v-2(n(x) \cdot v)n(x)) + (1-c)C_M \mathcal{M}(v) \int_{n(x) \cdot u > 0} f(t, x, u)(n(x) \cdot u) ~du.
\end{equation*} 

\noindent Hence 

\begin{equation*}
(\mathcal{M}^{\frac{1}{2}} h)|_{\gamma_{-}} = c\mathcal{M}^{\frac{1}{2}}(v)h(t,x,v-2(n(x) \cdot v)n(x)) + (1-c)C_M \mathcal{M}(v) \int_{n(x) \cdot u>0} \mathcal{M}^{\frac{1}{2}}(u) h(t,x,u)(n(x) \cdot u) ~du ,
\end{equation*}

\noindent which gives rise to 

\begin{equation*}
    \begin{aligned}
        h(t, x, v)|_{\gamma_{-}} =& ch(t,x,v-2(n(x) \cdot v)n(x)) + (1-c)C_M \mathcal{M}^{\frac{1}{2}}(v) \int_{n(x) \cdot u>0} \mathcal{M}^{\frac{1}{2}}(u) h(t,x,u) (n(x) \cdot u) ~du \\
        =& ch(t,x,v-2(n(x) \cdot v)n(x)) + (1-c)P_{\gamma}h(t, x, v).
    \end{aligned}
\end{equation*}

\noindent Then we have 
\begin{equation*}
    \begin{aligned}
        \int_{n(x) \cdot v<0} P_{\gamma}^2h(t,x,v) (n(x) \cdot v) ~dv =& -\int_{n(x) \cdot v>0} P_{\gamma}^2h(t,x,v) (n(x) \cdot v) ~dv
        = -\| P_{\gamma} h \|_{\gamma_{+}}^2,
    \end{aligned}
\end{equation*}

\noindent where we use the fact that $P_{\gamma}h(t,x,v)=P_{\gamma}h(t,x,-v)$. Let $v_{*}=v-2(n(x) \cdot v)n(x)$. Clearly the absolute value of the Jacobian of this transformation is 1, and $n(x)\cdot v =-n(x) \cdot v_{*}$. Hence we get 

\begin{equation*}
    \int_{n(x) \cdot v<0} h^2(t,x,v-2(n(x) \cdot v)n(x)) (n(x) \cdot v) ~dv = -\int_{n(x) \cdot v_{*}>0} h^2(t,x,v_{*}) (n(x) \cdot v_{*}) ~dv_{*} = -\| h\|^2_{\gamma_{+}},
\end{equation*}

\noindent and 

\begin{equation*}
    \begin{aligned}
        &\int_{n(x) \cdot v<0} h(t,x,v-2(n(x) \cdot v)n(x))P_{\gamma}h(t,x,v) (n(x) \cdot v) ~dv \\
        =& -\int_{n(x) \cdot v_{*}>0} h(t,x,v_{*})P_{\gamma}h(t,x,v_{*}) (n(x) \cdot v_{*}) ~dv_{*} \\
        =& -\int_{n(x) \cdot v_{*}>0} h(t,x,v_{*})C_M\mathcal{M}^{\frac{1}{2}}_{*}(n(x) \cdot v_{*}) ~dv_{*} \int_{n(x)\cdot u >0}\mathcal{M}^{\frac{1}{2}}(u) h(t,x,u) (n(x) \cdot u) ~du \\
        =& - \left(\int_{n(x) \cdot v>0} h(t,x,v)C_M^{\frac{1}{2}}\mathcal{M}^{\frac{1}{2}}(n(x) \cdot v) ~dv\right)^2 
        = - \| P_{\gamma}h\|^2_{\gamma_{+}},
    \end{aligned}
\end{equation*}

\noindent where we use the fact that $P_{\gamma}h(t,x,v)=P_{\gamma}h(t,x,v_{*})$. Collecting all these equalities, we have 

\begin{equation*}
    \int_{n(x) \cdot v<0} h^2(n(x) \cdot v) ~dv = -c^2\| h\|^2_{\gamma_{+}} -(1-c^2)\| P_{\gamma}h\|^2_{\gamma_{+}}.
\end{equation*}

\noindent And hence

\begin{equation*}
    \int_{\mathbb{R}^3} h^2 (n(x) \cdot v) ~dv = \| h\|^2_{\gamma_{+}} - c^2\| h\|^2_{\gamma_{+}} -(1-c^2)\| P_{\gamma}h\|^2_{\gamma_{+}} = (1-c^2)\| (I-P_{\gamma})h\|^2_{\gamma_{+}},
\end{equation*}

\noindent which leads to 

\begin{equation*}
    \langle v \cdot \nabla_x f, f \rangle = \frac{1}{2}\int_{\mathbb{R}^3} \int_{\partial \Omega} f^2 \mathcal{M}^{-1}(v)(v \cdot n(x)) ~d S_x dv = \frac{1-c^2}{2}|(I-P_{\gamma}) h|^{2}_{2,+} \geq 0.
\end{equation*} 

Now go back to \eqref{entropy equation}, we have 
\begin{equation*}
        \frac{1}{2} \frac{d}{dt} \| f \|^2 = - \frac{1-c^2}{2}|(I-P_{\gamma}) h |_{2,+}^2 + \langle \mathcal{L}(f), f \rangle 
        \leq  - \frac{1-c^2}{2}|(I-P_{\gamma}) h |_{2,+}^2 -\lambda \| f^{\perp} \|^2 ,
\end{equation*}

\noindent by Lemma \ref{microscopic coercivity}. This completes the proof. 

\qed

Now, we give the proof of Theorem \ref{exponential decay of f base case}. 

\bigskip

\textit{Proof of Theorem \ref{exponential decay of f base case}.} Since $H(f)=\frac{1}{2}\| f \|^2 + \eta\langle j, \nabla_x \phi \rangle_{d x} $, we have 

\begin{equation*}
    \begin{aligned}
     \frac{d}{dt} H(f)=& \frac{1}{2} \frac{d}{dt}\| f \|^2 + \eta \frac{d}{dt}\langle j, \nabla_x \phi \rangle_{dx} \\
     =& \frac{1}{2} \frac{d}{dt}\| f \|^2 + \eta \langle \partial_t j, \nabla_x \phi \rangle_{dx} + \eta \langle j, \partial_t \nabla_x \phi \rangle_{dx} \\
     =& \frac{1}{2} \frac{d}{dt}\| f \|^2 + \eta \langle -\nabla_x \cdot S - \nabla_x \rho + j^{\mathcal{L}}, \nabla_x \phi \rangle_{d x} + \eta\langle j, \partial_t \nabla_x \phi \rangle_{dx}.
    \end{aligned}
\end{equation*}

Let $T_1 := \frac{1}{2} \frac{d}{d t}\|f\|^2, ~T_2 := -\left\langle\nabla_x \cdot S, \nabla_x \phi\right\rangle_{d x}, ~T_3 :=-\left\langle\nabla_x \rho, \nabla_x \phi\right\rangle_{d x}, ~T_4 := \langle j^{\mathcal{L}} , \nabla_x \phi \rangle_{dx},  ~T_5:= \langle j, \partial_t \nabla_x \phi \rangle_{dx}$. Then $\frac{d}{dt} H(f) = T_1 +\eta(T_2 + T_3 + T_4 + T_5)$. We will do the following estimate for each $T_i$.  

For $T_1$, by Lemma \ref{entropy estimate}, we have 

\begin{equation}\label{estimate for T1}
    T_1 \leq - \frac{1-c^2}{2}| (I-P_{\gamma}) h|_{2,+}^{2} -\lambda \|f^{\perp} \|^2.
\end{equation}

For $T_3$, using integration by parts and the Neumann boundary condition on $\phi$, we have

\begin{equation}\label{estimate for T3}
    T_3 = -\langle \nabla_x \rho , \nabla_x \phi \rangle_{dx} = \langle \rho , \Delta \phi \rangle_{dx} +\int_{\partial \Omega} \rho (\nabla_x \phi \cdot n(x)) ~d S_x = -\| \rho \|_{dx}^2.
\end{equation}

For $T_4$, using Lemma \ref{estimate for j_L}, Lemma \ref{estimates_for_phi}, Cauchy-Schwarz and Young's inequality, we have 

\begin{equation}\label{estimate for T4}
    T_4 = \langle j^{\mathcal{L}}, \nabla_{x} \phi \rangle_{dx} \leq \| j^{\mathcal{L}}\|_{dx} \cdot \| \nabla_{x} \phi \|_{dx} \leq C_{\mathcal{L}}\| f^{\perp} \| \cdot C_{p}\| \rho \|_{dx} \leq \frac{C_{\mathcal{L}}C_{p}}{4\eta_3} \| f^{\perp} \|^2 + \eta_{3}C_{\mathcal{L}}C_{p} \| \rho \|^2_{dx},
\end{equation}

\noindent where $\eta_3>0$ is a free variable to be chosen later. 

For $T_5$, using Cauchy-Schwarz, Lemma \ref{estimates_for_phi} and Lemma \ref{moment estimates}, we have 
\begin{equation}\label{estimate for T5}
    T_5 = \langle j, \partial_t \nabla_x \phi \rangle_{dx} \leq \| j \|_{d x} \cdot \|\partial_t \nabla_x \phi\|_{dx} \leq \| j \|_{dx}^2 \leq 3\| f^{\perp} \|^2.
\end{equation}

Finally, for $T_2$, applying integration by parts leads to

\begin{equation*}
    T_2 = -\langle \nabla_x \cdot S, \nabla_x \phi \rangle_{dx} =-\sum_{i=1}^3 \langle \nabla_x \cdot S^i , \phi_{x_i} \rangle_{dx} =\sum_{i=1}^3 \langle S^i, \nabla_x \phi_{x_i} \rangle_{dx} - \sum_{i=1}^3 \int_{\partial \Omega} \phi_{x_i} (S^i \cdot n(x)) ~dS_x. 
\end{equation*}

\noindent For the term $\sum_{i=1}^3 \langle S^i, \nabla_x \phi_{x_i} \rangle_{dx}$, using Cauchy-Schwarz and Young's inequality, we have 

\begin{equation}\label{estimate for T_2 normal part}
    \begin{aligned}
        \sum_{i=1}^3 \langle S^i, \nabla_x \phi_{x_i} \rangle_{dx} \leq& \frac{1}{4 \eta_1} \sum_{i=1}^3 \| S^i \|_{d x}^2+\eta_1 \sum_{i=1}^3 \|\nabla_x \phi_{x_i} \|_{dx}^2 \\
        =& \frac{1}{4 \eta_1}\| S \|_{d x}^2 + \eta_1 \|\nabla_x^2 \phi \|_{d x}^2 \\
        \leq& \frac{D}{4 \eta_1} \|f^{\perp} \|^2+\eta_1 K^2 \| \rho \|_{dx}^2 ,\\
    \end{aligned}
\end{equation}

\noindent where $\eta_1>0$ is a free variable to be chosen later and $\nabla^2_{x} \phi = (\phi_{x_ix_j})$ is the Hessian matrix of $\phi$. Notice that we use Lemma \ref{moment estimates} and a regularity estimate for Poisson's equation in the last inequality above.

For the surface integral $\sum_{i=1}^3 \int_{\partial \Omega} \phi_{x_i} (S^i \cdot n(x)) ~dS_x = \int_{\partial \Omega} \sum_{i=1}^3 \phi_{x_i} (S^i \cdot n(x)) ~dS_x$, fix a point $x \in \partial\Omega$, we have the following
\begin{equation*}
    \begin{aligned}
        &\sum_{i=1}^3 \phi_{x_i}(S^i \cdot n(x)) 
        = S n(x) \cdot \nabla_x \phi \\
        =& (\tilde{S}n(x) - \rho n(x)) \cdot \nabla_x \phi \hspace{1cm} (\text{by the definitions of } S \text{ and } \tilde{S})\\
        =& \tilde{S}n(x) \cdot \nabla_x \phi  \hspace{1cm} (\text{by Neumann boundary condition on } \phi)\\
        =& \int_{\mathbb{R}^3}f(v \otimes v)n(x) \cdot \nabla_{x}\phi ~dv 
        = \int_{\mathbb{R}^3} h\mathcal{M}^{\frac{1}{2}}(n(x) \cdot v)(v \cdot \nabla_{x}\phi) ~dv \\
        =& \int_{n(x) \cdot v >0} h\mathcal{M}^{\frac{1}{2}}(n(x) \cdot v)(v \cdot \nabla_{x}\phi) ~dv + \int_{n(x) \cdot v <0} h\mathcal{M}^{\frac{1}{2}}(n(x) \cdot v)(v \cdot \nabla_{x}\phi) ~dv\\
        =& \int_{n(x) \cdot v >0} h\mathcal{M}^{\frac{1}{2}}(n(x) \cdot v)(v \cdot \nabla_{x}\phi) ~dv + c\int_{n(x) \cdot v <0} h(t,x,v-2(n(x) \cdot v)n(x))\mathcal{M}^{\frac{1}{2}}(n(x) \cdot v)(v \cdot \nabla_{x}\phi) ~dv \\
        &+ (1-c)\int_{n(x) \cdot v <0} P_{\gamma}h(t, x, v)|_{\gamma_{-}}\mathcal{M}^{\frac{1}{2}}(n(x) \cdot v)(v \cdot \nabla_{x}\phi) ~dv \\
        =& (1-c)\int_{n(x) \cdot v >0} h\mathcal{M}^{\frac{1}{2}}(n(x) \cdot v)(v \cdot \nabla_{x}\phi) ~dv + (1-c)\int_{n(x) \cdot v <0} P_{\gamma}h(t, x, v)|_{\gamma_{-}}\mathcal{M}^{\frac{1}{2}}(n(x) \cdot v)(v \cdot \nabla_{x}\phi) ~dv,
    \end{aligned}
\end{equation*}

\noindent where $f=\mathcal{M}^{\frac{1}{2}}h$ is defined as in the proof of Lemma \ref{entropy estimate}, and we use the fact that $v \cdot \nabla_{x}\phi = v_{*} \cdot \nabla_{x}\phi$ for $v_{*}= v-2(n(x) \cdot v)n(x)$ in the last equality. Hence the specular reflection part of the boundary vanishes and without loss of generality, we may assume $h$ satisfies the pure diffusive boundary condition. We may further assume(with a possible linear transformation) that $n(x)=(0,0,1)^{T}$. By the Neumann boundary condition, we must have $\phi_{x_{3}}=0$. Then write $h=P_{\gamma}h + (I-P_{\gamma})h$, we get 

\begin{equation*}
    \int_{\mathbb{R}^3} h\mathcal{M}^{\frac{1}{2}}(n(x) \cdot v)(v \cdot \nabla_{x}\phi) ~dv = \int_{\mathbb{R}^3} (P_{\gamma}h + (I-P_{\gamma})h)\mathcal{M}^{\frac{1}{2}}v_3(v_1\phi_{x_{1}}+v_{2}\phi_{x_{2}}) ~dv.
\end{equation*}

\noindent Notice that $P_{\gamma}h = g(x)\mathcal{M}^{\frac{1}{2}}$ for some function $g$ in $x$ by \eqref{projection on boundary}. So $\int_{\mathbb{R}^3} (P_{\gamma}h)\mathcal{M}^{\frac{1}{2}}v_3(v_1\phi_{x_{1}}+v_{2}\phi_{x_{2}}) ~dv=0$ by the fact that $\int_{\mathbb{R}^3} v_iv_j\mathcal{M}(v) ~dv=0 ~~~\forall i\neq j$. For the term $\int_{\mathbb{R}^3} ((I-P_{\gamma})h)\mathcal{M}^{\frac{1}{2}}v_3(v_1\phi_{x_{1}}+v_{2}\phi_{x_{2}}) ~dv$, first note that $h|_{\gamma_{-}}=P_{\gamma}h$ as shown in the proof of Lemma \ref{entropy estimate}. So $(I-P_{\gamma})h|_{\gamma_{-}}=0$ and hence 
\begin{equation*}
    \begin{aligned}
        \int_{\mathbb{R}^3} ((I-P_{\gamma})h)\mathcal{M}^{\frac{1}{2}}v_3(v_1\phi_{x_{1}}+v_{2}\phi_{x_{2}}) ~dv =& \int_{v_3>0} ((I-P_{\gamma})h)\mathcal{M}^{\frac{1}{2}}v_3(v_1\phi_{x_{1}}+v_{2}\phi_{x_{2}}) ~dv \\
        =& \int_{n(x) \cdot v >0} ((I-P_{\gamma})h)\mathcal{M}^{\frac{1}{2}}(v_1\phi_{x_{1}}+v_{2}\phi_{x_{2}}) (n(x) \cdot v)~dv.
    \end{aligned}
\end{equation*}

\noindent Using Cauchy-Schwarz with respect to $\langle \cdot , \cdot\rangle_{\gamma_{+}}$, we can derive 
\begin{equation*}
    |\int_{n(x) \cdot v >0} ((I-P_{\gamma})h)\mathcal{M}^{\frac{1}{2}}v_1\phi_{x_{1}} (n(x) \cdot v)~dv| = |\phi_{x_{1}}| \cdot | \langle (I-P_{\gamma})h, v_1\mathcal{M}^{\frac{1}{2}}\rangle_{\gamma_{+}} | \leq |\phi_{x_{1}}| \cdot \| (I-P_{\gamma})h \|_{\gamma_{+}} \cdot \| v_1\mathcal{M}^{\frac{1}{2}} \|_{\gamma_{+}}.
\end{equation*}

\noindent Clearly similar estimate holds for $|\int_{n(x) \cdot v >0} ((I-P_{\gamma})h)\mathcal{M}^{\frac{1}{2}}v_2\phi_{x_{2}} (n(x) \cdot v)~dv|$. We set 
\begin{equation*}
    C_{\gamma} = \| v_1\mathcal{M}^{\frac{1}{2}} \|_{\gamma_{+}} = \| v_2\mathcal{M}^{\frac{1}{2}} \|_{\gamma_{+}} = \int_{v_{3}>0} v_2^2v_3 \mathcal{M} ~dv < \infty.
\end{equation*}   

\noindent Collecting all these estimates, we have 
\begin{equation}\label{estimate for T_2 boundary part}
    \begin{aligned}
        -& \sum_{i=1}^3 \int_{\partial \Omega} \phi_{x_i} (S^i \cdot n(x)) ~dS_x \\
        =& - \int_{\partial \Omega} \int_{\mathbb{R}^3} h\mathcal{M}^{\frac{1}{2}}(n(x) \cdot v)(v \cdot \nabla_{x}\phi) ~dvdS_{x} \\
        \leq& \int_{\partial \Omega}  \int_{\mathbb{R}^3} |h\mathcal{M}^{\frac{1}{2}}(n(x) \cdot v)(v \cdot \nabla_{x}\phi) ~dv | dS_{x} \\
        \leq& (1-c)\int_{\partial \Omega} C_{\gamma}(|\phi_{x_1}|+|\phi_{x_2}|) \| (I-P_{\gamma}) h\|_{\gamma_{+}} ~dS_{x} \\
        \leq& 2(1-c)C_{\gamma} \| \nabla_{x}\phi\|_{L^2({\partial \Omega})} \cdot |(I-P_{\gamma}) h|_{2,+} \hspace{0.7cm} (\text{by Cauchy-Schwarz with respect to } \langle \cdot,\cdot \rangle_{L^2(\partial \Omega)}) \\
        \leq& 2(1-c)C_{\gamma}D_{\gamma}\| \nabla_{x}\phi \|_{H^{1}(\Omega)} \cdot |(I-P_{\gamma}) h|_{2,+} \hspace{0.7cm} (\text{by trace theorem}) \\
        \leq& 2(1-c)C_{\gamma}D_{\gamma}K\| \rho \|_{dx} \cdot |(I-P_{\gamma}) h|_{2,+} \hspace{0.7cm} (\text{by regularity estimate for Poisson's equation}) \\
        \leq& 2(1-c)C_{\gamma}D_{\gamma} K\eta_2\| \rho \|_{dx}^{2} + \frac{(1-c)C_{\gamma}D_{\gamma} K}{2\eta_2} |(I-P_{\gamma}) h|^{2}_{2,+} ,\hspace{0.7cm} (\text{by Young's inequality})
    \end{aligned}
\end{equation}

\noindent where $D_{\gamma}$ is the constant appearing in the trace theorem with respect to $\nabla_{x} \phi$, and $\eta_2>0$ is a free variable to be chosen. We therefore derive the following estimate for $T_2$
\begin{equation}\label{estimate for T2}
    T_{2} \leq \frac{D}{4 \eta_1} \|f^{\perp} \|^2+(\eta_1 K^2+2(1-c)C_{\gamma}D_{\gamma} K\eta_2) \| \rho \|_{dx}^2 + \frac{(1-c)C_{\gamma}D_{\gamma} K}{2\eta_2} | (I-P_{\gamma}) h|_{2,+}^{2}.
\end{equation}

Since $\eta_1, \eta_2, \eta_3$ are free, we set them to be equal so that $\eta_1 = \eta_2 =\eta_3= \tilde{\eta} >0$. Then collecting all the estimates \eqref{estimate for T1} - \eqref{estimate for T2}, we have 
\begin{equation*}
    \begin{aligned}
        \frac{d}{dt}H(f) \leq& \underbrace{\left(-\lambda +\eta(\frac{D}{4\tilde{\eta}}+ \frac{C_{\mathcal{L}}C_{p}}{4\tilde{\eta}} + 3)\right)}_{:= \alpha_{\eta,\tilde{\eta}}}\| f^{\perp} \|^2 + \underbrace{\eta\left(\tilde{\eta}K^2 + 2(1-c)\tilde{\eta}C_{\gamma}D_{\gamma}K  +\tilde{\eta}C_{\mathcal{L}}C_{p}-1\right)}_{:=\beta_{\eta,\tilde{\eta}}} \| \rho\|^2_{dx} \\
        &+ \underbrace{\left(-\frac{1-c^2}{2} + \frac{\eta (1-c)C_{\gamma}D_{\gamma}K}{2\tilde{\eta}}\right)}_{:= \delta_{\eta,\tilde{\eta}}} | (I-P_{\gamma}) h|_{2,+}^{2}.
    \end{aligned}
\end{equation*}

\noindent We first choose $\tilde{\eta}< \frac{1}{K^2 + 2(1-c)C_{\gamma}D_{\gamma}K + C_{\mathcal{L}}C_p}$ so that $\beta_{\eta,\tilde{\eta}}<0$ regardless of the choice of $\eta$. Once $\tilde{\eta}$ is fixed, we choose $\eta < \min(\frac{\lambda}{\frac{D}{4\tilde{\eta}}+ \frac{C_{\mathcal{L}}C_{p}}{4\tilde{\eta}} + 3}~, \frac{\tilde{\eta}(1+c)}{C_{\gamma}D_{\gamma}K})$ so that $\alpha_{\eta,\tilde{\eta}}<0, \delta_{\eta,\tilde{\eta}} \leq 0$. By Lemma \ref{equivalent norm}, $H(f)$ and $\| f \|^2$ are equivalent if $\eta < \frac{1}{2\sqrt{3}C_{p}}$. So we choose $ \eta < \min( \frac{\lambda}{\frac{D}{4\tilde{\eta}}+ \frac{C_{\mathcal{L}}C_{p}}{4\tilde{\eta}} + 3}, \frac{\tilde{\eta}(1+c)}{C_{\gamma}D_{\gamma}K}, \frac{1}{2\sqrt{3}C_{p}} )$, and let $\omega = -\max ( \alpha_{\eta,\tilde{\eta}}, \beta_{\eta,\tilde{\eta}} )>0$. Then 

\begin{equation*}
    \frac{d}{dt}H(f) \leq -\omega(\| f^{\perp} \|^2 + \| \rho \|^2_{dx}) = -\omega\| f \|^2 \leq -\frac{\omega}{C_{\eta}}H(f).
\end{equation*}

\noindent The result follows by using Gronwall's inequality and Lemma \ref{equivalent norm}.

\qed 
\section{generalization with potential} 

In this section, we will generalize the result of the previous section to include electrical potential in our model equation, which is just \eqref{linear semiconductor Boltzmann equation} and we restate it here: 

\begin{equation}\label{linear semiconductor with potential}
    \begin{cases}
         \partial_{t}f + v\cdot \nabla_{x}f + \nabla_{x}V \cdot \nabla_{v}f= \mathcal{L}(f) \\
         \ f(0,x,v) = f_{0}(x,v), \hspace{0.7cm} (x,v) \in \Omega \times \mathbb{R}^3.
    \end{cases}
\end{equation}

For simplicity, we assume $V$ is smooth and is invariant under time $t$, ie: $V=V(x) \in C^2(\Omega)$. Notice that adding any constant to $V$ will not change the electrical field $\nabla_{x}V$, hence we can normalize $V$ so that $\int_{\Omega} e^{V} ~dx=1$. This normalization makes the equilibrium state of the system unique. 

The main result of this section is an analogue of Theorem \ref{exponential decay of f base case} when electrical potential is involved. However, notice that when there is electrical potential, $\mathcal{M}(v)$ is no longer an equilibrium state to \eqref{linear semiconductor with potential}, since the left-hand side of \eqref{linear semiconductor with potential} is not $0$ when substituting $f=\mathcal{M}$. Instead, it can be checked that the equilibrium state is given by $F(x,v) = e^{V}\mathcal{M}$. So it is more natural to consider the norm $\| \cdot \|_{dud\nu}$ with respect to the weighted measure:
\begin{align*}
    &dud\nu:=F^{-1}dxdv = e^{-V}\mathcal{M}^{-1}dxdv, 
\end{align*}rather than the norm $\| \cdot\|_{dxd\nu}$ as in the previous section. Note that since $V \in C^{2}(\Omega)$ is continuous over a bounded domain, we have $c_{V} \leq e^{-V} \leq C_{V}$ for some $0 < c_{V} \leq C_{V}$. Hence the norms $\| \cdot \|_{dud\nu}$ and $\| \cdot \|_{dxd\nu}$ are equivalent, with
\begin{equation}\label{equivalence between dudnu and dxdnu}
    c_{V}\| \cdot \|_{dxd\nu}^2 \leq \| \cdot \|_{dud\nu}^2 \leq C_{V}\| \cdot \|_{dxd\nu}^2.
\end{equation}

\noindent Although the norms are equivalent, the proofs in this section are easier using $\| \cdot \|_{dud\nu}$. So we let $\langle \cdot,\cdot \rangle:=\langle \cdot,\cdot \rangle_{dud\nu}$, $\| \cdot \|:=\| \cdot \|_{dud\nu}$(note that these are different from the notations we used in the last section). We now state the main result of this section. 

\begin{theorem}\label{exponential decay of f with potential}
    Let $f$ be a solution to \eqref{linear semiconductor with potential}, where the cross-section $\sigma(v,v_{*})$ of $\mathcal{L}$ satisfies \textbf{Assumption 1} and \textbf{Assumption 2}, and the electrical potential $V=V(x) \in C^2(\Omega)$. Suppose $f$ is given the Maxwell boundary condition on $\partial \Omega \times \mathbb{R}^3$ with parameter $c$. Then there exists $C>0$ and $\tau >0$ such that 

    \begin{equation}\label{exponential decay of f equation with potential}
        \| f(t) - \mu_{f}e^{V}\mathcal{M}\| \leq C \| f_{0} -\mu_{f}e^{V}\mathcal{M} \| ~e^{-\tau t},
    \end{equation}

\noindent $\forall t \geq 0$, where $\mu_{f} = \iint_{\Omega \times \mathbb{R}^3} f_{0} ~dxdv$.  
\end{theorem}

To prove Theorem \ref{exponential decay of f with potential}, we follow similar methodology as the proof of Theorem \ref{exponential decay of f base case}. Thanks to the equivalence of norms, we can apply most of the lemmas in the previous section directly and the proofs are essentially the same. The only exception is the proof of Lemma \ref{entropy estimate}, where we need a more careful treatment of the terms involving potentials. So we present a new proof of the lemma as follows

\begin{lemma}\label{entropy estimate with potential}
    \textnormal{(Entropy estimate with potential).} Let $f$ be a solution to \eqref{linear semiconductor with potential} with Maxwell boundary condition, and let $\lambda>0$ be defined as in \textbf{Assumption 1}. Then 
    \begin{equation}\label{entropy estimate equation with potential}
        \frac{1}{2} \frac{d}{dt} \| f \|^2 \leq - \frac{1-c^2}{2}| (I-P_{\gamma}) h|_{2,+,du}^{2} -\lambda \| f^{\perp} \|^2 ,
    \end{equation}

    \noindent where $h:=\mathcal{M}^{-\frac{1}{2}}f$ and $|h|_{2,+,du}^{2} := \int_{\partial \Omega}\int_{n(x)\cdot v>0} |h(t,x,v)|^2 e^{-V} (n(x)\cdot v)~dvd S_{x}$. 
\end{lemma}

\textit{Proof}. Take $\langle, f \rangle$ on both sides of \eqref{linear semiconductor with potential}, we get:
\begin{equation}\label{entropy equation with potential}
    \frac{1}{2} \frac{d}{dt}\|f\|^2 + \langle v \cdot \nabla_x f, f \rangle + \langle \nabla_{x}V \cdot \nabla_{v}f, f \rangle = \langle \mathcal{L}(f), f \rangle .
\end{equation}

For the term $\langle v \cdot \nabla_x f, f \rangle$, using integration by parts with respect to $\nabla_{x}$, we have
\begin{equation*}
    \begin{aligned}
        \langle v \cdot \nabla_x f, f \rangle = & \int_{\mathbb{R}^3} \int_{\Omega}(\nabla_x \cdot (vf)) f e^{-V}\mathcal{M}^{-1} ~dxdv \\
        = & - \int_{\mathbb{R}^3} \int_{\Omega} (vf) \cdot \nabla_x (fe^{-V}) \mathcal{M}^{-1} ~dxdv + \int_{\mathbb{R}^3} \int_{\partial \Omega} f^2 e^{-V} \mathcal{M}^{-1}(n(x) \cdot v) ~dS_{x}dv \\
        = & - \underbrace{\int_{\mathbb{R}^3} \int_{\Omega} (vf) \cdot (e^{-V}\nabla_x f) \mathcal{M}^{-1} ~dxdv}_{=\langle v \cdot \nabla_{x}f , f \rangle} - \int_{\mathbb{R}^3} \int_{\Omega} (vf) \cdot (-fe^{-V}\nabla_{x}V) \mathcal{M}^{-1} ~dxdv \\
        & + \int_{\mathbb{R}^3} \int_{\partial \Omega} f^2 e^{-V} \mathcal{M}^{-1}(n(x) \cdot v) ~dS_{x}dv.
\end{aligned}
\end{equation*}

For the term $\langle \nabla_{x}V \cdot \nabla_{v}f, f \rangle$, using integration by parts with respect to $\nabla_{v}$, we have

\begin{equation*}
    \begin{aligned}
        \langle \nabla_{x}V \cdot \nabla_{v}f, f \rangle = & \int_{\mathbb{R}^3} \int_{\Omega}(\nabla_{v}\cdot f \nabla_x V) f e^{-V} \mathcal{M}^{-1} ~dxdv \\
        = & - \int_{\mathbb{R}^3} \int_{\Omega} (f \nabla_{x}V) \cdot \nabla_v (f \mathcal{M}^{-1})e^{-V} ~dxdv \\
        = & - \underbrace{\int_{\mathbb{R}^3} \int_{\Omega} (f \nabla_{x}V) \cdot (\mathcal{M}^{-1}\nabla_{v}f)e^{-V} ~dxdv}_{=\langle \nabla_{x}V \cdot \nabla_{v}f, f \rangle}  - \int_{\mathbb{R}^3} \int_{\Omega} (f \nabla_{x}V) \cdot (f v\mathcal{M}^{-1})e^{-V} ~dxdv,
\end{aligned}
\end{equation*}

\noindent where we use the fact that $\nabla_{v}\mathcal{M}^{-1}= v\mathcal{M}^{-1}$ for the second term in the last equality above. Combining these equations, we get 

\begin{equation*}
    \langle v \cdot \nabla_x f, f \rangle + \langle \nabla_{x}V \cdot \nabla_{v}f, f \rangle = \frac{1}{2}\int_{\partial \Omega}e^{-V} \left(\int_{\mathbb{R}^3} f^2 \mathcal{M}^{-1}(n(x) \cdot v) ~dv \right) dS_{x}.
\end{equation*}

Let $f=\mathcal{M}^{\frac{1}{2}}h$ and follow exactly the same arguments as the proof of Lemma \ref{entropy estimate}, we can derive 

\begin{equation*}
    \int_{\mathbb{R}^3} f^2 \mathcal{M}^{-1}(n(x) \cdot v) ~dv = (1-c^2)\| (I-P_{\gamma})h\|^2_{\gamma_{+}}.
\end{equation*}

\noindent Then the result follows using Lemma \ref{microscopic coercivity}. 

\qed

We now present the proof of Theorem \ref{exponential decay of f with potential}. 

\textit{Proof of Theorem \ref{exponential decay of f with potential}.} Since $e^{V}\mathcal{M}$ is an equilibrium to \eqref{linear semiconductor with potential}, we can replace $f$ by $f+ae^{V}\mathcal{M}$ for any $a \in \mathbb{R}$ without changing the equation or the boundary condition. So without loss of generality, we can assume $\mu_{f}=0$. Then the proof follows the same structure as the proof of Theorem \ref{exponential decay of f base case}. First, we consider the same moments of $f$ as defined in \eqref{moments of f}. Then multiplying \eqref{linear semiconductor with potential} by $1,v$ and taking integration, we get the moment equations with potential 
\begin{equation}\label{moment equations with potential}
   \begin{cases}
        \partial_t \rho + \nabla_x \cdot j=0 \\
        \partial_t j+\nabla_x \cdot(S + \rho I) -(\nabla_{x}V) \rho = j^{\mathcal{L}},
    \end{cases}
\end{equation}

\noindent where we use the fact that $\int_{\mathbb{R}^3} \nabla_{v}f ~dv=0$ and $\int_{\mathbb{R}^3} (\nabla_{x}V \cdot \nabla_{v}f)v ~dv= -\int_{\mathbb{R}^3} (\nabla_{x}V) f ~dv = -(\nabla_{x}V)\rho$. 

We construct the ``same" modified entropy functional $H(f)=\frac{1}{2}\| f \|^2 + \eta \langle j, \nabla_{x}\phi \rangle_{du}$ as \eqref{def of H(f)}, where $du=e^{-V}dx$ and $\phi$ solves the same Poisson equation as \eqref{Poisson equation for phi} with the same Neumann boundary condition. The estimates for the moments and $\phi$ follow exactly the same as Lemma \ref{moment estimates} and \ref{estimates_for_phi}, where the proofs do not depend on the potential $V$. $H(f)$ and $\| f \|^2$ are again equivalent when $\eta>0$ is small enough, with 

\begin{equation*}
    c_{\eta}\| f \|^2 \leq H(f) \leq C_{\eta}\| f \|^2,
\end{equation*}

\noindent where $c_{\eta}=\frac{1}{2} - \frac{\sqrt{3}\eta C_{p}}{c_{V}}$ and $C_{\eta}= \frac{1}{2} + \frac{\sqrt{3}\eta C_{p}}{c_{V}}$. 

Consider $\frac{d}{dt}H(f) = \frac{1}{2}\frac{d}{dt}\| f \|^2 + \eta \frac{d}{dt} \langle j, \nabla_{x}\phi \rangle_{du}$, we have

\begin{equation*}
    \frac{d}{dt} H(f) = T_1 +\eta(T_2 + T_3 + T_4 + T_5 + T_6),
\end{equation*}

\noindent where $T_1 := \frac{1}{2} \frac{d}{d t}\|f\|^2, ~T_2 := -\left\langle\nabla_x \cdot S, \nabla_x \phi\right\rangle_{du}, ~T_3 :=-\left\langle\nabla_x \rho, \nabla_x \phi\right\rangle_{du}, ~T_4 := \langle j^{\mathcal{L}} , \nabla_x \phi \rangle_{du},  ~T_5:= \langle j, \partial_t \nabla_x \phi \rangle_{du}$, and $~T_6:= \langle (\nabla_{x}V) \rho, \nabla_{x} \phi \rangle_{du}$. 

The estimate for $T_1$ is given by Lemma \ref{entropy estimate with potential}, i.e. 

\begin{equation*}
    T_1 \leq - \frac{1-c^2}{2}| (I-P_{\gamma}) h|_{2,+,du}^{2} -\lambda \| f^{\perp} \|^2.
\end{equation*}

\noindent Consider $T_3 + T_6$, we have 

\begin{equation*}
    \begin{aligned}
         T_3 +T_6 =&  -\left\langle\nabla_x \rho, \nabla_x \phi\right\rangle_{du} + \langle (\nabla_{x}V) \rho, \nabla_{x} \phi \rangle_{du} \\
         =& -\int_{\Omega} \big(\nabla_{x}\rho \cdot \nabla_{x}\phi \big)e^{-V}~dx + \int_{\Omega} \big( \nabla_{x}V \cdot \nabla_{x}\phi \big)\rho e^{-V}~dx \\
         =& \int_{\Omega}\rho \big( -\rho e^{-V} -(\nabla_{x}\phi \cdot \nabla_{x}V)e^{-V}\big) ~dx + \int_{\Omega} \big( \nabla_{x}V \cdot \nabla_{x}\phi \big)\rho e^{-V}~dx \\
         =& -\| \rho\|^{2}_{du}.
    \end{aligned}
\end{equation*}

\noindent For the estimates of $T_4, T_5$, we apply the same arguments as \eqref{estimate for T4}, \eqref{estimate for T5} and the fact that $\| \cdot\|_{du}$ is equivalent to $\| \cdot \|_{dx}$. Then we get

\begin{equation*}
    \begin{cases}
        T_4 \leq \| j^{\mathcal{L}}\|_{du}\| \nabla_{x}\phi\|_{du} \leq C_{\mathcal{L}}C_{p}C_{V}\| f^{\perp}\|_{dxd\nu}\| \rho\|_{dx} \leq \frac{C_{\mathcal{L}}C_{p}\tilde{C}_{V}}{4\eta_{3}}\| f^{\perp}\|^2 +\eta_{3}C_{\mathcal{L}}C_{p}\tilde{C}_{V}\| \rho\|^2_{du}\\
        T_5 \leq \| j\|_{du}\| \partial_t \nabla_x \phi\|_{du} \leq 3C_{V}\| f^{\perp}\|^2_{dxd\nu} \leq 3\tilde{C}_{V}\| f^{\perp}\|^2,
    \end{cases}
\end{equation*}

\noindent where $\tilde{C}_{V}=\frac{C_{V}}{c_{V}}$. Finally for $T_2$, we have  
\begin{equation*}
    \begin{aligned}
        T_2 = -\langle \nabla_x \cdot S, \nabla_x \phi \rangle_{du} =&-\sum_{i=1}^3 \langle \nabla_x \cdot S^i , \phi_{x_i} \rangle_{du} \\
        =& -\sum_{i=1}^3\int_{\Omega}(\nabla_x \cdot S^i)\phi_{x_{i}}e^{-V}~dx \\
        =& \sum_{i=1}^3\int_{\Omega} S^{i}\cdot (\nabla_{x}\phi_{x_i}e^{-V}-\phi_{x_i}e^{-V}\nabla_{x}V) ~dx- \sum_{i=1}^3 \int_{\partial \Omega} \phi_{x_i}e^{-V} (S^i \cdot n(x)) ~dS_x\\
        =&\sum_{i=1}^3 \langle S^i, \nabla_x \phi_{x_i} \rangle_{du} -\sum_{i=1}^3 \langle S^i, \phi_{x_i}\nabla_{x}V \rangle_{du}  - \sum_{i=1}^3 \int_{\partial \Omega} \phi_{x_i}e^{-V} (S^i \cdot n(x)) ~dS_x.
    \end{aligned}
\end{equation*}

\noindent For the term $\sum_{i=1}^3 \langle S^i, \nabla_x \phi_{x_i} \rangle_{du} -\sum_{i=1}^3 \langle S^i, \phi_{x_i}\nabla_{x}V \rangle_{du}$, using similar arguments as \eqref{estimate for T_2 normal part} and assume $\| V\|_{C^{2}(\Omega)} \leq D_V$(use the fact that $V \in C^2(\Omega)$), we have 

\begin{equation*}
    \sum_{i=1}^3 \langle S^i, \nabla_x \phi_{x_i} \rangle_{du} -\sum_{i=1}^3 \langle S^i, \phi_{x_i}\nabla_{x}V \rangle_{du} \leq \frac{2D\tilde{C}_{V}}{4 \eta_1} \|f^{\perp} \|^2+\eta_1 K^2(1+D_V)\tilde{C}_{V} \| \rho \|_{du}^2.
\end{equation*}

\noindent To estimate the term $- \sum_{i=1}^3 \int_{\partial \Omega} \phi_{x_i}e^{-V} (S^i \cdot n(x)) ~dS_x$, we apply the same arguments as \eqref{estimate for T_2 boundary part}(replacing $dx$ by $du$) to derive 
\begin{equation*}
    - \sum_{i=1}^3 \int_{\partial \Omega} \phi_{x_i}e^{-V} (S^i \cdot n(x)) ~dS_x \leq 2\eta_{2}(1-c)C_{\gamma}D_{\gamma} K\tilde{C}_{V}\| \rho \|_{du}^{2} + \frac{(1-c)C_{\gamma}D_{\gamma} K\tilde{C}_{V}}{2\eta_2} |(I-P_{\gamma}) h|^{2}_{2,+,du}.
\end{equation*}

Collecting all these estimates and setting $\eta_1=\eta_2=\eta_3=\tilde{\eta}$, we have 

\begin{equation*}
    \begin{aligned}
        \frac{d}{dt} H(f) \leq& \left( -\lambda + \frac{2\eta D\tilde{C}_{V}}{4 \tilde{\eta}} + \frac{\eta C_{\mathcal{L}}C_{p}\tilde{C}_{V}}{4\tilde{\eta}} + 3\eta\tilde{C}_{V} \right) \| f^{\perp} \|^2 + \eta ( -1 +\tilde{\eta} K^2(1+D_V)\tilde{C}_{V} + 2\tilde{\eta}(1-c)C_{\gamma}D_{\gamma} K\tilde{C}_{V} 
        \\
        &+ \tilde{\eta}C_{\mathcal{L}}C_{p}\tilde{C}_{V} )\| \rho \|_{du}^2 + \left( - \frac{1-c^2}{2} + \frac{\eta(1-c)C_{\gamma}D_{\gamma} K\tilde{C}_{V}}{2\tilde{\eta}}\right)|(I-P_{\gamma}) h|^{2}_{2,+,du}.
    \end{aligned}
\end{equation*}

\noindent Then, we apply the same arguments as the last paragraph of the proof of Theorem \ref{exponential decay of f base case}: we set $\tilde{\eta}$ and $\eta$ to be small enough and the result follows using a Gronwall-type inequality. 

\qed

\section{generalization with uncertainty}

In this section, we will generalize the previous results to incorporate randomness into our model. We mainly consider randomness arising from the collision cross-section $\sigma$ and the initial data $f_{0}$, and we assume that all randomness can be characterized by a finite-dimensional random vector $z \in \mathbb{R}^{d}$ with known probability distribution. For more details on the parameterization of random inputs by a finite-dimensional random vector, we refer the reader to the appendix. To simplify our arguments, We assume $z \in I_{z} \subset \mathbb{R}$ is a single random variable(cases for higher dimensions can be generalized without any essential difficulty), and we also ignore the external potential in this section. Hence we consider the following equation 

\begin{equation}\label{linear semiconductor with UQ}
    \begin{cases}
         \partial_{t}f + v\cdot \nabla_{x}f = \mathcal{L}_{z}(f) \\
         \ f(0,x,v,z) = f_{0}(x,v,z), \hspace{0.7cm} (x,v,z) \in \Omega \times \mathbb{R}^3 \times I,
    \end{cases}
\end{equation}

\noindent with the same Maxwell boundary condition as in the previous sections. The subscript of $\mathcal{L}_{z}$ is used to denote the randomness of the collision operator, and we have $\mathcal{L}_{z}(f) := \int_{\mathbb{R}^3}\sigma(v,v_{*},z)(\mathcal{M}f_{*} - \mathcal{M}_{*}f) ~dv_{*}$. We assume both $\sigma$ and $f_0$ are sufficiently smooth in $z$, and we want to know how the smoothness is propagated through time. In particular, we need to verify the regularity of $f$ with respect to $z$ before we apply any numerical method to solve the stochastic system \eqref{linear semiconductor with UQ}, such as the gPC method discussed in \cite{HJ16}, \cite{LiuJin2018}. To this end, We let $g_{l}$ denote the lth derivative in $z$ of $f$: 

\begin{equation}\label{def of g_l}
    g^{l} = \frac{\partial f}{\partial z^{l}}
\end{equation}

\noindent and our goal of this section is to study the long time behavior of $\| g^{l} \|$ for $\| \cdot \|=\| \cdot\|_{dxd\nu}$ defined in section 3. To begin with, take $\frac{\partial}{\partial z^{l}}$ on both sides of \eqref{linear semiconductor with UQ} and apply the general Leibniz formula to yield

\begin{equation*}
    \begin{aligned}
        \partial_{t}g^{l} + v\cdot \nabla_{x}g^{l} =& \int_{\mathbb{R}^3} \frac{\partial}{\partial z^{l}} (\sigma(v,v_{*},z)(\mathcal{M}f_{*} - \mathcal{M}_{*}f)) ~dv_{*} \\
        =& \int_{\mathbb{R}^3} \sum_{k=0}^{l}\binom{l}{k} \frac{\partial\sigma}{\partial z^{l-k}} \cdot \frac{\partial(\mathcal{M}f_{*} - \mathcal{M}_{*}f)}{\partial z^{k}} ~dv_{*} \\
        =& \sum_{k=0}^{l}\binom{l}{k}\int_{\mathbb{R}^3}  \frac{\partial\sigma}{\partial z^{l-k}} (\mathcal{M}g^{k}_{*} - \mathcal{M}_{*}g^{k}) ~dv_{*} \\
        =& \mathcal{L}_{z}(g^l)+ \sum_{k=0}^{l-1} \binom{l}{k}\mathcal{L}_{z}^{l-k}(g^{k}),
    \end{aligned} 
\end{equation*}

\noindent where $\mathcal{L}_{z}^{k}(f) = \int_{\mathbb{R}^3}  \partial_{z}^{k}\sigma(v,v_{*},z)(\mathcal{M}f_{*} - \mathcal{M}_{*}f) ~dv_{*}$, $\partial_{z}^{k}\sigma(v,v_{*},z)=\frac{\partial\sigma(v,v_{*},z)}{\partial z^{k}}$. Set $\mathcal{S}^{l}=\sum_{k=0}^{l-1} \binom{l}{k}\mathcal{L}_{z}^{l-k}(g^{k})$, then $g^{l}$ satisfies a linear semiconductor Boltzmann equation with source term arising from lower derivatives of $f$, i.e.:

\begin{equation}\label{linear semiconductor equation with source}
    \begin{cases}
        \partial_{t}g^{l} + v\cdot \nabla_{x}g^{l} = \mathcal{L}_{z}(g^l) + \mathcal{S}^l \\
        \ g^{l}_{0}(x,v,z):= g^{l}(0,x,v,z) = \frac{\partial f_{0}}{\partial z^{l}}(x,v,z), \hspace{0.7cm} (x,v,z) \in \Omega \times \mathbb{R}^3 \times I.
    \end{cases}
\end{equation}

\noindent Note that $g^{l}$ satisfies the same Maxwell boundary condition as $f$, which can be easily checked by taking $\frac{\partial}{\partial z^{l}}$ on both sides of \eqref{maxwell_bc}. Using \eqref{local conservation of mass}, we have $\int_{\mathbb{R}^3}\mathcal{S}^{l} ~dv=0$. So the total mass $\iint_{\Omega \times \mathbb{R}^3} g^{l} ~dvdx$ is conserved pointwisely in $z$. Hence we can assume $\iint_{\Omega \times \mathbb{R}^3} g^{l} ~dvdx=0$ by similar arguments in the previous sections. We want to analyze the long time behavior of $g^{l}$, and our main result of this section is an analogue of Theorem \ref{exponential decay of f base case} for the derivatives of $f$ with respect to $z$, which is stated as follows:  

\begin{theorem}\label{exponential decay of g_l}
    Let $g^{l}$ be the $l$-th derivative of $f$ with respect to $z$, satisfying \eqref{linear semiconductor equation with source} Maxwell boundary condition with parameter $c$. Suppose \textbf{Assumption $1^{*},2^{*}, 3$} hold(given below), then there exists $a>0$ and a polynomial $G^{l}(t)$ in $t$ of degree $l$ with coefficients that are nonnegative linear combinations of $\| g^{k}(0) \|$ for $ 0 \leq k \leq l$, such that $\forall t \geq 0$, we have

    \begin{equation}\label{exponential decay of g_l equation}
        \| g^{l}(t) \| \leq ~e^{-at}G^{l}(t).
    \end{equation}
\end{theorem}

To prove Theorem \ref{exponential decay of g_l}, we can mostly follow the same arguments as in section 3, with two special considerations. First, with the additional random variable $z$, all the assumptions we made in section 3 need to be true pointwise in $z$, and we shall set uniform(in $z$) bounds to simplify our arguments. To be specific, we have the following stochastic version of assumptions 1-2, with an additional assumption on the uniform boundedness of the operators $\mathcal{L}_{z}^{l}$:

\bigskip
\noindent \textbf{Assumption $1^{*}$} 

\begin{equation}\label{H1*}
    \sigma(v,v_{*},z) \geq \lambda(z) \geq \lambda >0.
\end{equation}

\bigskip
\noindent \textbf{Assumption $2^{*}$}

\begin{equation}\label{H2*}
    \iint_{\mathbb{R}^3 \times \mathbb{R}^3} \partial_{z}^{l}\sigma(v,v_{*},z)^{2}\mathcal{M}\mathcal{M}_{*}(|v|^2 + |v_{*}|^2) ~dvdv_{*} \leq C_{l}(z) \leq C \ <\infty \hspace{0.8cm} \forall l \geq 0.
\end{equation}

\bigskip
\noindent \textbf{Assumption 3}

\begin{equation}\label{H5}
    \| \mathcal{L}_{z}^{l} \| = \tilde{C}_{l}(z) \leq \tilde{C}<\infty \hspace{0.8cm} \forall l \geq 0,
\end{equation}

\bigskip
\noindent where $\lambda, C,\tilde{C}>0$ are constants independent of $z$.

Another consideration is that we need to analyze how the source term $\mathcal{S}^{l}$ affects our arguments in the previous sections. In particular, our entropy estimate shall be adjusted in the presence of a source $\mathcal{S}^{l}$: 

\begin{lemma}\label{entropy estimate with source}
    \textnormal{(Entropy estimate with source).} Let $g^{l}$ be a solution to \eqref{linear semiconductor equation with source} with Maxwell boundary condition, and let $\lambda>0$ be defined as in \textbf{Assumption $1^{*}$}. Let $h^{l}:=\mathcal{M}^{-\frac{1}{2}}g^{l}$, then 

    \begin{equation}\label{entropy estimate equation with source}
        \frac{1}{2} \frac{d}{dt} \| g^{l} \|^2 \leq - \frac{1-c^2}{2}| (I-P_{\gamma}) h^{l}|_{2,+}^{2}-\lambda \| g^{l,\perp} \|^2 + \| \mathcal{S}^{l} \| \| g^{l} \|.
    \end{equation}
\end{lemma}

\textit{Proof.} Take $\langle, g^{l} \rangle$ on both sides of \eqref{linear semiconductor equation with source}, we get:

\begin{equation}\label{entropy equation with source}
    \frac{1}{2} \frac{d}{dt}\|g^{l}\|^2 + \langle v \cdot \nabla_x g^{l}, g^{l} \rangle = \langle \mathcal{L}_{z}(g^{l}), g^{l} \rangle + \langle \mathcal{S}^{l} , g^{l} \rangle.
\end{equation}

Using exactly the same arguments as in the proof of Lemma \ref{entropy estimate}, we have 
$\langle v \cdot \nabla_x g^{l}, g^{l} \rangle = \frac{1-c^2}{2}| (I-P_{\gamma}) h^{l}|_{2,+}^{2} $. Using \textbf{Assumption $1^{*}$} and Lemma \ref{microscopic coercivity}, we have $\langle \mathcal{L}_{z}(g^{l}), g^{l} \rangle \leq -\lambda\| g^{l,\perp} \|^{2}$. Then using Cauchy-Schwarz for the term $\langle \mathcal{S}^{l}, g^{l} \rangle$ yields the result.  

\qed

Since the source term $\mathcal{S}^{l}$ depends on the derivatives of lower orders, it is natural to suggest that there should be some recursive relations involved in deriving our final estimates. We therefore propose the following lemma, which will be useful later for the proof of our main result. 

\begin{lemma}\label{recursion estimates}
    Let $\{ h^{l}(t)\}_{l \geq 0}$ be a sequence of nonnegative $C^{1}$ functions in $t$ such that 
    \begin{equation*}
        \begin{cases}
            \frac{d}{dt}h^{0} \leq -a h^{0}    \\
            \frac{d}{dt}h^{l} \leq -a h^{l} +\sum_{k=0}^{l-1}b_{k}^{l}h^{k} ,\hspace{0.8cm} \forall l \geq 1,
        \end{cases}
    \end{equation*}

    \noindent where $a>0, b_{k}^{l} \geq 0 ~\forall l \geq 1, 0 \leq k <l$. Then we have 

    \begin{equation*}
        h^{l}(t) \leq e^{-at}G^{l}(t),
    \end{equation*}

    \noindent where $G^{l}(t)$ is a polynomial in $t$ of degree $l$ with coefficients that are linear combinations of $h^{k}(0)$ for $0 \leq k \leq l$, and $G^{l}(t)$ is given by the recursive formula 
    \begin{equation}\label{recursion estimates formula}
        \begin{cases}
            G^{0}(t) = h^{0}(0) \\
            G^{l}(t) = h^{l}(0) + \sum_{k=0}^{l-1}b^{l}_{k}\int_{0}^{t}G^{k}(s) ~ds, \hspace{0.8cm} \forall l \geq 1.
        \end{cases}
    \end{equation}
\end{lemma}

\textit{Proof.} Consider $\frac{d}{dt}(e^{at}h^{l})$, we have 
\begin{equation*}
        \frac{d}{dt}(e^{at}h^{l}) = a e^{at}h^{l} + e^{at}\frac{d}{dt}h^{l} \leq a e^{at}h^{l} +e^{at}\left(-a h^{l} +\sum_{k=0}^{l-1}b_{k}^{l}h^{k}\right) = e^{at}\sum_{k=0}^{l-1}b_{k}^{l}h^{k}.
\end{equation*}

\noindent Integrate both sides of the inequality above from $0$ to $t$ gives 
\begin{equation*}
    e^{at}h^{l}(t) - h^{l}(0) \leq \int_{0}^{t} e^{as}\sum_{k=0}^{l-1}b_{k}^{l}h^{k}(s) ~ds.
\end{equation*}

\noindent Then clearly we have $h^{0}(t) \leq e^{-at}h^{0}(0)=e^{-at}G^{0}(t)$. For $h^{1}(t)$, we can derive the following 
\begin{align*}
    &h^{1}(t) \leq e^{-at}\left(h^{1}(0) + \int_{0}^{t} e^{as} \cdot b^{1}_{0}h^{0}(s) ~ds\right) \leq e^{-at}\left(h^{1}(0) + \int_{0}^{t} e^{as} \cdot b^{1}_{0}e^{-as}h^{0}(0) ~ds\right) \\
    &\leq e^{-at}\left(h^{1}(0) + b^{1}_{0}h^{0}(0)t\right)=e^{-at}G^{1}(t).
\end{align*}

\noindent Then the result follows by a simple induction on $l$. 

\qed

We now present the proof of Theorem \ref{exponential decay of g_l}.  

\textit{Proof of Theorem \ref{exponential decay of g_l}.} With the source term, the moment equations now become 
\begin{equation}\label{moment equations with source}
    \begin{cases}
        \partial_t \rho + \nabla_x \cdot j=0 \\
        \partial_t j+\nabla_x \cdot(S + \rho I) = j^{\mathcal{L}} + j^{\mathcal{S}^{l}},
    \end{cases}
\end{equation}

\noindent where $j^{\mathcal{S}^{l}} = \int_{\mathbb{R}^3} v\mathcal{S}^{l} ~dv$, and $\rho, j , S$ are now macroscopic quantities associated to $g^{l}$. We again construct the modified entropy functional $H(g^{l}) = \frac{1}{2} \| g^{l}\|^2 + \eta \langle j,\nabla_{x}\phi \rangle_{dx}$ as in \eqref{def of H(f)}, where $\phi$ satisfies the same Poisson equation \eqref{Poisson equation for phi} with the same Neumann boundary condition. The same moment estimates(Lemma \ref{moment estimates}) and estimates for $\phi$(Lemma \ref{estimates_for_phi}) apply pointwisely in $z$, since the proofs do not depend on the source term $\mathcal{S}^{l}$. $H(g^{l})$ is again equivalent to $\| g^{l}\|^2_{dxd\nu}$ when $\eta< \frac{1}{2\sqrt{3}C_{p}}$, where $C_{p}$ is the constant for Poincare's inequality, which is independent of $z$. Then we have the same constants of equivalence $c_{\eta}$ and $C_{\eta}$ as in \eqref{equivalent norm equation}, both independent of $z$. Consider the time evolution of $H(g^{l})$, we get 

\begin{equation*}
    \frac{d}{dt}H(g^{l}) = T_1 +\eta(T_2 + T_3 +T_4 +T_5 +T_6),
\end{equation*}

\noindent where $T_1 - T_5$ are defined in the same way as in the proof of Theorem \ref{exponential decay of f base case}, with estimates for $T_2 - T_5$ given exactly the same as \eqref{estimate for T3}-\eqref{estimate for T2}. The estimate for $T_1$ is given by Lemma \ref{entropy estimate with source} and $T_6 := \langle j^{\mathcal{S}^{l}},\nabla_{x}\phi \rangle_{dx}$. Using Lemma \ref{estimates_for_phi} and Cauchy-Schwarz, we have 

\begin{equation}\label{estimate for T6 with randomness}
    T_6 = \langle j^{\mathcal{S}^{l}},\nabla_{x}\phi \rangle_{dx} \leq C_{p}\| j^{\mathcal{S}^{l}} \|_{dx}\| \rho\|_{dx} \leq C_{p}\| j^{\mathcal{S}^{l}} \|_{dx}\| g^{l} \|.
\end{equation}

Then collecting all the estimates, we derive 

\begin{equation*}
    \begin{aligned}
        \frac{d}{dt}H(f) \leq& \left(-\lambda +\eta(\frac{D}{4\tilde{\eta}}+ \frac{C_{\mathcal{L}}C_{p}}{4\tilde{\eta}} + 3)\right)\| g^{l,\perp} \|^2 + \eta\left(\tilde{\eta}K^2 + 2(1-c)\tilde{\eta}C_{\gamma}D_{\gamma}K  +\tilde{\eta}C_{\mathcal{L}}C_{p}-1\right)\| \rho\|^2_{dx} \\
        &+ \left(-\frac{1-c^2}{2} + \frac{\eta (1-c)C_{\gamma}D_{\gamma}K}{2\tilde{\eta}}\right)| (I-P_{\gamma}) h^l|_{2,+}^{2}+ \| \mathcal{S}^{l}\|\| g^{l}\| + C_{p}\| j^{\mathcal{S}^{l}} \|_{dx}\| g^{l} \|.
    \end{aligned}
\end{equation*}

\noindent where we set $\eta_1=\eta_2=\eta_3=\tilde{\eta}$. Following the same arguments as the proof of Theorem \ref{exponential decay of f base case}, we can properly choose $\tilde{\eta}$ and $\eta$ so that for some $\omega>0$,
\begin{equation*}
    \frac{d}{dt}H(g^{l}) \leq -\omega \| g^{l}\|^2 + \| \mathcal{S}^{l}\|\| g^{l}\| + C_{p}\| j^{\mathcal{S}^{l}} \|_{dx}\| g^{l} \|.
\end{equation*}

\noindent  Notice that all the parameters at this stage are independent of the random variable $z$. Set $\varphi^{l}=\sqrt{H(g^{l})}$, Recall that when $\eta$ is small enough, $H(g^{l})$ is equivalent to $\| g^{l} \|^2$(hence $\varphi^{l}$ is equivalent to $\| g^{l}\|$) with constants of equivalence $c_{\eta}$ and $C_{\eta}$($\sqrt{c_{\eta}}$ and $\sqrt{C_{\eta}}$). Using the fact that $\frac{d}{dt}H(g^{l})=2\varphi^{l}\cdot \frac{d}{dt}\varphi^{l}$, we have 

\begin{equation*}
    \frac{d}{dt}\varphi^{l} \leq -\frac{\omega}{2C_{\eta}}\varphi^{l} + \frac{1}{2\sqrt{c_{\eta}}}\| \mathcal{S}^{l} \| + \frac{C_{p}}{2\sqrt{c_{\eta}}}\| j^{\mathcal{S}^{l}}\|_{dx}.
\end{equation*}

\noindent Since $\mathcal{S}^{l}=\sum_{k=0}^{l-1} \binom{l}{k}\mathcal{L}_{z}^{l-k}(g^{k})$, using \textbf{Assumption 3}, we get 

\begin{equation*}
    \| \mathcal{S}^{l} \| = \| \sum_{k=0}^{l-1} \binom{l}{k}\mathcal{L}_{z}^{l-k}(g^{k}) \| \leq \sum_{k=0}^{l-1} \binom{l}{k} \| \mathcal{L}_{z}^{l-k}(g^{k})\| \leq \sum_{k=0}^{l-1} \binom{l}{k}\tilde{C} \| g^{k}\|.
\end{equation*}

\noindent Since $j^{\mathcal{S}^{l}} = \int_{\mathbb{R}^3} v\mathcal{S}^{l} ~dv = \sum_{k=0}^{l-1} \binom{l}{k} \int_{\mathbb{R}^3} v\mathcal{L}_{z}^{l-k}(g^{k}) ~dv$, using \textbf{Assumption $2^{*}$} and Lemma \ref{estimate for j_L}, we have

\begin{equation*}
    \| j^{\mathcal{S}^{l}} \|_{dx} \leq \sum_{k=0}^{l-1} \binom{l}{k} \| \int_{\mathbb{R}^3} v\mathcal{L}_{z}^{l-k}(g^{k}) ~dv \|_{dx} \leq \sum_{k=0}^{l-1} \binom{l}{k} C_{\mathcal{L}}\| g^{k,\perp} \| \leq \sum_{k=0}^{l-1} \binom{l}{k} C_{\mathcal{L}}\| g^{k} \|.
\end{equation*}

\noindent Collecting all these estimates, we get 

\begin{equation*}
    \frac{d}{dt}\varphi^{l} \leq -\frac{\omega}{2C_{\eta}}\varphi^{l} + \frac{1}{2\sqrt{c_{\eta}}}\sum_{k=0}^{l-1} \binom{l}{k}\tilde{C} \| g^{k}\| + \frac{C_{p}}{2\sqrt{c_{\eta}}} \sum_{k=0}^{l-1} \binom{l}{k} C_{\mathcal{L}}\| g^{k} \| \leq -\frac{\omega}{2C_{\eta}}\varphi^{l} + \sum_{k=0}^{l-1} \binom{l}{k}\frac{\tilde{C}+C_{p}C_{\mathcal{L}}}{2c_{\eta}} \varphi^{k}.
\end{equation*}

Set $a=\frac{\omega}{2C_{\eta}}$ and $b^{l}_{k} = \binom{l}{k}\frac{\tilde{C}+C_{p}C_{\mathcal{L}}}{2c_{\eta}}$, then the inequality above is reduced to the form presented in Lemma \ref{recursion estimates}. The result follows by applying Lemma \ref{recursion estimates} and the equivalence between $\varphi^{l}$ and $\| g^{l}\|$. 

\qed

\section{Conclusion} 

In this paper, we establish hypocoercivity for the semiconductor Boltzmann equation under the Maxwell boundary condition by constructing a modified entropy Lyapunov functional. This functional is proved to be equivalent to some weighted norm of the corresponding function space, and it can be shown to dissipate along the solutions. Then the exponential decay to the equilibrium state of the system follows by a Gronwall-type inequality. We also generalize our arguments to situations where uncertainties in our model are presented, and the hypocoercivity method we have established is adopted to analyze the regularity of the solutions along the random space. In the future work, we will further incorporate uncertainty quantification into our model by enlarging sources of uncertainties, e.g. uncertainty arising from wall temperature. 
\appendix

\section{Parametrization of random inputs}

In this appendix, we discuss how to manage the randomness existing in our model. The key step is to parameterize the random inputs by a finite set of independent random variables to make computational simulations plausible. If the random inputs are already given in the form of finitely many random parameters with a proper probability distribution, e.g., jointly Gaussian, then parametrization is straightforward, e.g., using Cholesky decomposition. However, in many cases, the random inputs are formulated by random processes, which are often characterized by a continuous index $t \in T$. Then we need to apply dimension reduction techniques to approximate the processes using finitely many random variables. One of the most widely used techniques in this regard is the Karhunen-Loeve(KL) expansion, see \cite{Hajek}, \cite{Xiu}, \cite{Sullivan}. 

For a random process $ \{Y_{t}(\omega) \}_{t \in T}$(we use the notation $Y_t$ as an abbreviation) with mean $\mu_{Y}(t)$ and autocorrelation function $R_{Y}(t,s)=\mathbb{E}[Y_tY_s]$, its KL expansion, if exists, is given by 

\begin{equation}\label{def of KL expansion}
    Y_{t}(\omega) = \sum_{i=1}^{\infty} \psi_{i}(t) Y_{i}(\omega),
\end{equation}

\noindent where the series converges in the mean square sense(m.s.), $\{\psi_{i}\}$ is an orthonormal family of functions in $L^2(T)$ and $Y_{i}$'s are mutually orthogonal, i.e., $\mathbb{E}[Y_{i}Y_{j}]=0 ~~\forall i \neq j$. The existence of KL expansion for $Y_{t}$ is guaranteed by Mercer's theorem(see \cite{Hajek} for more details), provided that the random process $Y_{t}$ is m.s. continuous, i.e., $R_Y(t,s)$ is continuous over $T \times T$. 

We now analyze how to derive $\psi_{i}$ and $Y_i$ if the KL expansion for $Y_t$ exists. The KL expansion for $Y_t$ can be viewed as an analogue of decomposing $f \in L^2(T)$ with respect to an orthonormal family $\{\psi_{i}\}$, i.e., 

\begin{equation*}
    f=\sum_{i=1}^{\infty} f_{i}\psi_{i},
\end{equation*}  

\noindent where $f_i = \int_{T} f\psi_{i} ~dt \in \mathbb{R}$ are the Fourier coefficients. In the KL setting, once an orthonormal family $\{\psi_{i}\}$ is chosen, we can define the "Fourier coefficients" $Y_i$ analogously by 

\begin{equation}\label{Fourier coefficient for Y_i}
    Y_i(\omega) = \int_{T} Y_{t}(\omega)\psi_{i}(t) ~dt.
\end{equation}

\noindent The only difference is that the coefficients for the KL expansion are random variables instead of real numbers. The challenge now is to pick $\{\psi_{i}\}$ properly so that the coefficients $\{Y_{i}\}$ are mutually orthogonal. This can be accomplished by the following lemma, whose proof is included in \cite{Hajek}.  

\begin{lemma}\label{eigenvalue problem for psi_i}
    Suppose $Y_t$ is m.s. continuous and \eqref{def of KL expansion} holds for $Y_t$ with $\{\psi_i\}$ orthonormal and $\{Y_i\}$ not necessarily mutually orthogonal. Then it is a KL expansion(i.e., $Y_i$'s are mutually orthogonal) if and only if $\psi_i$'s are eigenfunctions of $R_Y$: 

    \begin{equation*}
        R_Y(\psi_i) = \lambda_i \psi_i, 
    \end{equation*}

    \noindent where $R_Y$ is an operator on $L^2(T)$ given by $R_Y(\psi)(t)=\int_{T} R_Y(t,s)\psi(s) ~ds$ for any $\psi(t) \in L^2(T)$. In case \eqref{def of KL expansion} is a KL expansion, the eigenvalues are given by $\lambda_i = \mathbb{E}[|Y_i|^2]$. 
\end{lemma}

\noindent In summary, if $Y_t$ is m.s. continuous, its KL expansion can be established by first solving an eigenvalue problem related to the autocorrelation function $R_Y(t,s)$ to obtain an orthonormal family $\{ \psi_i \}$, followed by a computation of the "Fourier coefficients" $Y_i(\omega)$ associated to this orthonormal family. Once the KL expansion is established, the analysis of $Y_t$ can be naturally transformed into the analysis of the coefficients $Y_i$. 

For practical purposes, we need to truncate the series appearing in \eqref{def of KL expansion} to obtain a finite dimensional parametrization of the random process, i.e., 

\begin{equation}\label{KL expansion with truncation}
    Y_{t}(\omega) \approx \sum_{i=1}^{d} \psi_{i}(t) Y_{i}(\omega), \hspace{0.8cm} d \geq 1.
\end{equation}

\noindent In most situations, the eigenvalues $\lambda_{i}$ as appeared in lemma \ref{eigenvalue problem for psi_i} will decay as $i$ increases. Hence we can choose the truncation order $d$ based on the decay rate of the eigenvalues. For more details, we refer the readers to \cite{Xiu}. Once \eqref{KL expansion with truncation} is established, we can represent the random process $Y_t$ by finitely many orthogonal random variables $Y_i$ as we desire. Note that in general $Y_i$'s are not mutually independent, unless additional assumptions on $Y_t$ are made, e.g., $Y_t$ is a Gaussian process with zero mean. We will not pursue further in this direction and shall be content with finite representation of $Y_t$ by orthogonal random variables. Some remarks are in order. 

\begin{remark}
    The condition of m.s. continuity to guarantee the existence of a KL expansion is not very restrictive. Many random processes we use for modeling, e.g., Brownian motion and Poisson process, satisfy this property. 
\end{remark}

\begin{remark}
    We assume the distribution of the random process $Y_t$ is prescribed. Hence, we can derive the probability distributions of $Y_{i}$ using \eqref{Fourier coefficient for Y_i}. This is especially straightforward and useful when $Y_{t}$ is a Gaussian process. 
\end{remark}

\begin{remark}
    The truncated KL expansion \eqref{KL expansion with truncation} identifies the "most accurate" $d$-dimensional approximation of $Y_t$ in the sense that it minimizes $\mathbb{E}[\| Y_t - Z_t\|^2]$ over all $d$-dimensional random processes $Z_t$. A random process $Z_t$ is said to be $d$-dimensional if it has the form $Z_{t}(\omega) = \sum_{i=1}^{d} \phi_{i}(t) Z_{i}(\omega)$ for any $d$ random variables $Z_1,...,Z_d$ and functions $\phi_1,...,\phi_d$. 
\end{remark}

\section{Proofs of preliminary lemmas}

\noindent \textit{proof of lemma \ref{microscopic coercivity}}. By definition, we have

\begin{equation*}
    \begin{aligned}
        \langle \mathcal{L}(f) , f \rangle_{d\nu} =& \int \int_{\mathbb{R}^3 \times \mathbb{R}^3} f\mathcal{M}^{-1} \sigma(v,v_{*})(\mathcal{M}f_{*} - \mathcal{M}_{*}f) ~dvdv_{*} \\
        =& \int \int_{\mathbb{R}^3 \times \mathbb{R}^3} \frac{f}{\mathcal{M}} \sigma(v,v_{*})\mathcal{M}\mathcal{M}_{*}\left(\frac{f_{*}}{\mathcal{M}_{*}} - \frac{f}{\mathcal{M}}\right) ~dvdv_{*} \\
        =& \int \int_{\mathbb{R}^3 \times \mathbb{R}^3} \frac{f_{*}}{\mathcal{M}_{*}} \sigma(v,v_{*})\mathcal{M}\mathcal{M}_{*}\left(\frac{f}{\mathcal{M}} - \frac{f_{*}}{\mathcal{M}_{*}}\right) ~dvdv_{*},
    \end{aligned}
\end{equation*}

\noindent where we exchange $v$ and $v_{*}$ to get the last equality above. Combining the second and third equality above, we have 

\begin{equation}\label{H theorem for L}
    \langle \mathcal{L}(f) , f \rangle_{d\nu} = -\frac{1}{2} \int \int_{\mathbb{R}^3 \times \mathbb{R}^3} \sigma(v,v_{*})\mathcal{M}\mathcal{M}_{*}\left(\frac{f}{\mathcal{M}} - \frac{f_{*}}{\mathcal{M}_{*}}\right)^2 ~dvdv_{*} \leq 0 ,
\end{equation}

\noindent which is the `H-theorem' for $\mathcal{L}$. Notice that \eqref{H theorem for L} holds for any $f \in L^2_{d\nu}$. In particular, it holds for $f^{\perp}$. Then, using Assumption 1, we have 

\begin{equation*}
    \begin{aligned}
        \langle \mathcal{L}(f^{\perp}) , f^{\perp} \rangle_{d\nu} =& -\frac{1}{2} \int \int_{\mathbb{R}^3 \times \mathbb{R}^3} \sigma(v,v_{*})\mathcal{M}\mathcal{M}_{*}\left(\frac{f^{\perp}}{\mathcal{M}} - \frac{f_{*}^{\perp}}{\mathcal{M}_{*}}\right)^2 ~dvdv_{*} \\
        \leq& -\frac{\lambda}{2} \int \int_{\mathbb{R}^3 \times \mathbb{R}^3}\mathcal{M}\mathcal{M}_{*}\left(\frac{f^{\perp}}{\mathcal{M}} - \frac{f_{*}^{\perp}}{\mathcal{M}_{*}}\right)^2 ~dvdv_{*} \\
        =& -\frac{\lambda}{2} \int \int_{\mathbb{R}^3 \times \mathbb{R}^3} \left(\frac{\mathcal{M}_{*}(f^{\perp})^2}{\mathcal{M}} -2f^{\perp}f^{\perp}_{*} + \frac{\mathcal{M}(f_{*}^{\perp})^2}{\mathcal{M}_{*}}\right) ~dvdv_{*} \\
        =& -\lambda \| f^{\perp}\|^2_{d\nu},
    \end{aligned}
\end{equation*}

\noindent where we use the fact that $\int_{\mathbb{R}^3} f^{\perp} ~dv=0$ and $\int_{\mathbb{R}^3} \mathcal{M}(v) ~dv=1$. Then notice that 

\begin{equation*}
    \langle \mathcal{L}(f) , f \rangle_{d\nu} = \langle \mathcal{L}(\rho \mathcal{M} + f^{\perp}) , \rho \mathcal{M} + f^{\perp} \rangle_{d\nu} = \langle \mathcal{L}(f^{\perp}) , f^{\perp} \rangle_{d\nu},
\end{equation*}

\noindent by \eqref{kernel of L} and \eqref{local conservation of mass}. This completes the proof.

\qed 

\bigskip

\noindent \textit{proof of lemma \ref{estimate for j_L}}. For each fixed $ x \in \Omega $, we use the notation $j^{\mathcal{L}} = j^{\mathcal{L}}(x)$ for simplicity. Then 

\begin{equation*}
    |j^{\mathcal{L}}| = |\int_{\mathbb{R}^3} v\mathcal{L}(f) ~dv | \leq \int_{\mathbb{R}^3} |v| |\mathcal{L}(f)| ~dv
\end{equation*}

\noindent and 

\begin{equation*}
    \begin{aligned}
        |\mathcal{L}(f)| = |\int_{\mathbb{R}^3}\sigma(v,v_{*})(\mathcal{M}f_{*} - \mathcal{M}_{*}f) ~dv_{*}| \leq& \int_{\mathbb{R}^3}\sigma(v,v_{*})(\mathcal{M}|f_{*}| + \mathcal{M}_{*}|f|) ~dv_{*} \\
        =& \int_{\mathbb{R}^3}\sigma(v,v_{*})\mathcal{M}|f_{*}| ~dv_{*} + \int_{\mathbb{R}^3}\sigma(v,v_{*})\mathcal{M}_{*}|f| ~dv_{*}.
    \end{aligned}
\end{equation*}

Consider the term $\int_{\mathbb{R}^3} |v| (\int_{\mathbb{R}^3}\sigma(v,v_{*})\mathcal{M}|f_{*}| ~dv_{*})~dv$, we have 

\begin{equation*}
    \begin{aligned}
        \int_{\mathbb{R}^3} |v| \left(\int_{\mathbb{R}^3}\sigma(v,v_{*})\mathcal{M}|f_{*}| ~dv_{*}\right)~dv =& \int_{\mathbb{R}^3} |v| \langle \sigma(v,v_{*})\mathcal{M}\mathcal{M_{*}}, |f_{*}| \rangle_{d\nu_{*}}~dv \\
        \leq& \int_{\mathbb{R}^3} |v| \| \sigma(v,v_{*})\mathcal{M}\mathcal{M_{*}} \|_{d\nu_{*}} \|f_{*} \|_{d\nu_{*}} ~dv \hspace{1 cm} \textnormal{by Cauchy-Schwarz} \\
        =& \| f \|_{d\nu} \langle |v| \| \sigma(v,v_{*})\mathcal{M}\mathcal{M_{*}} \|_{d\nu_{*}},\mathcal{M} \rangle_{d\nu} \\
        \leq & \| f \|_{d\nu} \|  |v| \| \sigma(v,v_{*})\mathcal{M}\mathcal{M_{*}} \|_{d\nu_{*}} \|_{d\nu} \hspace{1 cm} \textnormal{by Cauchy-Schwarz} \\
        =& \| f \|_{d\nu} \cdot \left(\int\int_{\mathbb{R}^3 \times \mathbb{R}^3} |v|^2 \sigma(v,v_{*})^{2}\mathcal{M}\mathcal{M}_{*} ~dv_{*}dv\right)^{\frac{1}{2}}.
    \end{aligned}
\end{equation*}

\noindent Similarly, we have 

\begin{equation*}
     \int_{\mathbb{R}^3} |v| \left(\int_{\mathbb{R}^3}\sigma(v,v_{*})\mathcal{M}_{*}|f| ~dv_{*}\right)~dv \leq \| f \|_{d\nu} \cdot \left(\int\int_{\mathbb{R}^3 \times \mathbb{R}^3} |v_{*}|^2 \sigma(v,v_{*})^{2}\mathcal{M}\mathcal{M}_{*} ~dv_{*}dv\right)^{\frac{1}{2}},
\end{equation*}

\noindent and hence 

\begin{equation*}
    |j^{\mathcal{L}}| \leq \| f \|_{d\nu} \cdot \left((\int\int_{\mathbb{R}^3 \times \mathbb{R}^3} |v|^2 \sigma(v,v_{*})^{2}\mathcal{M}\mathcal{M}_{*} ~dv_{*}dv)^{\frac{1}{2}} + (\int\int_{\mathbb{R}^3 \times \mathbb{R}^3} |v_{*}|^2 \sigma(v,v_{*})^{2}\mathcal{M}\mathcal{M}_{*} ~dv_{*}dv)^{\frac{1}{2}}\right).
\end{equation*}

Then the result follows using Assumption 2 and the fact that $\sqrt{a}+\sqrt{b} \leq \sqrt{2(a+b)}$. 

\qed 

\bigskip

\noindent \textit{proof of lemma \ref{null-flux property}}. By \eqref{maxwell_bc}, we have

\begin{equation*}
    \begin{aligned}
    \int_{\mathbb{R}^3} f(t, x, v) (n(x) \cdot v) ~dv = & \int_{n(x) \cdot v>0} f(t, x, v)(n(x) \cdot v) ~dv + \int_{n(x) \cdot v <0} f(t, x, v)(n(x) \cdot v) ~dv \\ 
    =& \int_{n(x) \cdot v>0} f(t, x, v)(n(x) \cdot v) ~dv + c\int_{n(x) \cdot v<0} f(t,x,v_{*})(n(x)\cdot v)~dv\\
    & + (1-c)\int_{n(x) \cdot v<0} C_{\mu} \mathcal{M}(v) \int_{n(x) \cdot u > 0} f(t, x, u)(n(x) \cdot u) ~du \cdot(n(x) \cdot v) ~dv \\
    =& \int_{n(x) \cdot v>0} f(t, x, v)(n(x) \cdot v) ~dv - c\int_{n(x) \cdot v_{*}>0} f(t,x,v_{*})(n(x)\cdot v_{*})~dv_{*} \\
    & -(1-c)\int_{n(x) \cdot u>0} f(t, x, u)(n(x) \cdot u) ~du \cdot \int_{n(x) \cdot v<0} C_{\mu} \mathcal{M}(v)|n(x) \cdot v| ~dv \\
    =& 0,
    \end{aligned}
\end{equation*}

\noindent where $v_{*}=v-2(n(x)\cdot v)n(x)$, and we use \eqref{unity of c_mu} to derive the last equality. This completes the proof. 

\qed

\bibliographystyle{spmpsci}
\bibliography{Ref.bib}

\begin{thebibliography}{10}
\providecommand{\url}[1]{{#1}}
\providecommand{\urlprefix}{URL }
\expandafter\ifx\csname urlstyle\endcsname\relax
  \providecommand{\doi}[1]{DOI~\discretionary{}{}{}#1}\else
  \providecommand{\doi}{DOI~\discretionary{}{}{}\begingroup \urlstyle{rm}\Url}\fi

\bibitem{A15}
Alonso, R.: Boltzmann-type equations and their applications.
\newblock Publica\c{c}\~{o}es Matem\'{a}ticas do IMPA. [IMPA Mathematical Publications]. Instituto Nacional de Matem\'{a}tica Pura e Aplicada (IMPA), Rio de Janeiro (2015).
\newblock 30${\sp{{}}{\rm{o}}}$ Col\'{o}quio Brasileiro de Matem\'{a}tica. [30th Brazilian Mathematics Colloquium]

\bibitem{BCMT2023}
Bernou, A., Carrapatoso, K., Mischler, S., Tristani, I.: Hypocoercivity for kinetic linear equations in bounded domains with general maxwell boundary condition.
\newblock Anal. Non Linéaire \textbf{40}, 287--338 (2023)

\bibitem{BHR2020}
Bessemoulin-Chatard, M., Herda, M., Rey, T.: Hypocoercivity and diffusion limit of a finite volume scheme for linear kinetic equations.
\newblock Mathematics of Computation \textbf{89}(323), 1093--1133 (2020)

\bibitem{briant2015}
Briant, M.: From the boltzmann equation to the incompressible navier-stokes equations on the torus: A quantitative error estimate.
\newblock Journal of Differential Equations \textbf{259}, 6072--6141 (2015)

\bibitem{briant2016asymptotic}
Briant, M., Guo, Y.: Asymptotic stability of the boltzmann equation with maxwell boundary conditions.
\newblock Journal of Differential Equations \textbf{261}(12), 7000--7079 (2016)

\bibitem{chen}
Chen, H.: {Cercignani-Lampis boundary in the Boltzmann theory}.
\newblock Kinetic \& Related Models \textbf{13}(3), 549--597 (2020)

\bibitem{chen2025boltzmann}
Chen, H., Duan, R.: Boltzmann equation with mixed boundary condition.
\newblock SIAM Journal on Mathematical Analysis \textbf{57}(3), 3297--3334 (2025)

\bibitem{CKP2025}
Chen, H., Klingenberg, C., Pirner, M.: {BGK model for rarefied gas in a bounded domain}.
\newblock arXiv preprint arXiv:2502.03096  (2025)

\bibitem{LiuDaus2019}
Daus, E.S., Jin, S., Liu, L.: Spectral convergence of the stochastic galerkin approximation to the boltzmann equation with multiple scales and large random perturbation in the collision kernel.
\newblock Kinetic and Related Models \textbf{12}, 909--922 (2019)

\bibitem{DM2015}
Dolbeault, J., Mouhot, C., Schmeiser, C.: Hypocoercivity for linear kinetic equations conserving mass.
\newblock Trans. Amer. Math. Soc. \textbf{367}(6), 3807--3828 (2015)

\bibitem{esposito2013non}
Esposito, R., Guo, Y., Kim, C., Marra, R.: Non-isothermal boundary in the boltzmann theory and fourier law.
\newblock Communications in Mathematical Physics \textbf{323}(1), 177--239 (2013)

\bibitem{esposito2018stationary}
Esposito, R., Guo, Y., Kim, C., Marra, R.: Stationary solutions to the boltzmann equation in the hydrodynamic limit.
\newblock Annals of PDE \textbf{4}(1), 1 (2018)

\bibitem{Guo2010}
Guo, Y.: Decay and continuity of the boltzmann equation in bounded domains.
\newblock Archive for rational mechanics and analysis \textbf{197}, 713--809 (2010)

\bibitem{guo2020landau}
Guo, Y., Hwang, H.J., Jang, J.W., Ouyang, Z.: The landau equation with the specular reflection boundary condition.
\newblock Archive for Rational Mechanics and Analysis \textbf{236}(3), 1389--1454 (2020)

\bibitem{Hajek}
Hajek, B.: Random Processes for Engineers.
\newblock Cambridge University Press (2015)

\bibitem{herau2006}
H{\'e}rau, F.: Hypocoercivity and exponential time decay for the linear inhomogeneous relaxation boltzmann equation.
\newblock Asymptot. Anal. \textbf{46}, 349--359 (2006)

\bibitem{herau2018}
H{\'e}rau, F.: Introduction to hypocoercive methods and applications for simple linear inhomogeneous kinetic models.
\newblock In: Lectures on the Analysis of Nonlinear Partial Differential Equations. Part 5, \emph{Morningside Lect. Math.}, vol.~5, pp. 119--147 (2018)

\bibitem{HN2004}
H{\'e}rau, F., Nier, F.: Isotropic hypoellipticity and trend to equilibrium for the fokker-planck equation with a high-degree potential.
\newblock Arch. Ration. Mech. Anal. \textbf{171}, 151--218 (2004)

\bibitem{HJ16}
Hu, J., Jin, S.: A stochastic {G}alerkin method for the {B}oltzmann equation with uncertainty.
\newblock J. Comput. Phys. \textbf{315}, 150--168 (2016)

\bibitem{JinPareschi}
Jin, S., Pareschi, L. (eds.): Uncertainty quantification for hyperbolic and kinetic equations, \emph{SEMA-SIMAI Springer Series}, vol.~14.
\newblock Springer (2017)

\bibitem{Jungel}
Jungel, A.: Transport Equations for Semiconductors.
\newblock Springer-Verlag, Berlin Heidelberg (2009)

\bibitem{LW2017}
Li, Q., Wang, L.: Uniform regularity for linear kinetic equations with random input based on hypocoercivity.
\newblock SIAM/ASA Journal on Uncertainty Quantification \textbf{5}(1), 1193--1219 (2017)

\bibitem{LiuJin2018}
Liu, L., Jin, S.: Hypocoercivity based sensitivity analysis and spectral convergence of the stochastic galerkin approximation to collisional kinetic equations with multiple scales and random inputs.
\newblock Multiscale Model. Simul. \textbf{16}, 1085--1114 (2018)

\bibitem{LQ2024-1}
Liu, L., Qi, K.: Convergence of the fourier-galerkin spectral method for the boltzmann equation with uncertainties.
\newblock Commun. Math. Sci. \textbf{22}(7), 1897--1925 (2024)

\bibitem{LQ2024-2}
Liu, L., Qi, K.: Spectral convergence of a semi-discretized numerical system for the spatially homogeneous boltzmann equation with uncertainties.
\newblock SIAM/ASA J. Uncertain. Quantif. \textbf{12}(3), 812--841 (2024)

\bibitem{mischler2010}
Mischler, S.: Kinetic equations with maxwell boundary conditions.
\newblock In: Annales scientifiques de l'Ecole normale sup{\'e}rieure, vol.~43, pp. 719--760 (2010)

\bibitem{MN2006}
Mouhot, C., Neumann, L.: Quantitative perturbative study of convergence to equilibrium for collisional kinetic models in the torus.
\newblock Nonlinearity \textbf{19}, 969--998 (2006)

\bibitem{parUQ}
Pareschi, L.: An introduction to uncertainty quantification for kinetic equations and related problems.
\newblock In: Trails in kinetic theory, \emph{SEMA SIMAI Springer Ser.}, vol.~25, pp. 141--181. Springer, Cham (2021)

\bibitem{Sullivan}
Sullivan, T.: Introduction to Uncertainty Quantification.
\newblock Springer, Cham (2015)

\bibitem{Villani02}
Villani, C.: A review of mathematical topics in collisional kinetic theory.
\newblock In: S.~Friedlander, D.~Serre (eds.) Handbook of Mathematical Fluid Mechanics, vol.~I, pp. 71--305. North-Holland (2002)

\bibitem{villani2009}
Villani, C.: Hypocoercivity, \emph{Mem. Amer. Math. Soc.}, vol. 202 (2009).
\newblock 184 pp.

\bibitem{WL2025}
Wan, J., Liu, L.: Error estimates of asymptotic-preserving neural networks in approximating stochastic linearized boltzmann equation.
\newblock arXiv preprint arXiv:2503.01643  (2025)

\bibitem{Xiu}
Xiu, D.: Numerical Methods for Stochastic Computations.
\newblock Princeton University Press, New Jersey (2010)

\end{thebibliography}

\end{document}